\documentclass[11pt]{amsart}

\usepackage{graphicx}
\usepackage[all]{xy}
\usepackage{amssymb,extarrows}

\newtheorem{thm}{Theorem}[section]
\newtheorem{cor}[thm]{Corollary}
\newtheorem{lem}[thm]{Lemma}
\newtheorem{prop}[thm]{Proposition}
\newtheorem{rem}[thm]{Remark}
\newtheorem{defn}[thm]{Definition}

\newcommand{\ra}{\longrightarrow}
\newcommand{\A}{\Lambda}
\newcommand{\s}{\blacksquare}
\newcommand{\xra}{\xlongrightarrow}
\newcommand{\End}{\operatorname{End}}
\newcommand{\add}{\operatorname{add}}
\newcommand{\Ima}{\operatorname{Im}}
\newcommand{\Cok}{\operatorname{Coker}}
\newcommand{\Ker}{\operatorname{Ker}}
\newcommand{\Hom}{\operatorname{Hom}}

\makeatletter
\begin{document}
\title{From CM-finite to CM-free}
\author{Fan Kong, Pu Zhang$^*$}
\thanks{$^*$ The corresponding author.}
\thanks{Supported by the NSF China (11271251), and Doc. Fund of Ministry of Education (20120073110058).} \maketitle
\begin{abstract} The aim of this paper is twofold. On one hand, we prove a slight
generalization of the stability for Gorenstein categories in [SWSW]
and [Huang]; and show that the relative Auslander algebra of a
CM-finite algebra is CM-free. On the other hand, we describe the
bounded derived category, and the Gorenstein defect category
introduced in [BJO], via Gorenstein-projective objects; and we show
that the Gorenstein defect category of a CM-finite algebra is
triangle-equivalent to the singularity category of its relative
Auslander algebra.
\end {abstract}



\section{\bf Introduction}

\subsection{} M. Auslander and M. Bridger [AB] introduced the modules of $G$-dimension zero over two-sided noetherian
rings. E. E. Enochs and O. M. G. Jenda [EJ1] generalized this
concept and introduced Gorenstein-projective modules over any ring.
This class of modules enjoys pleasant stable properties, becomes a
main ingredient in the relative homological algebra, and widely used
in the representation theory of algebras and algebraic geometry. See
e.g. [AM], [AR], [BGS], [Buch], [EJ2], [Hap2], [Hol], [Kn]. Recent
studies show that they are important not only to Gorenstein
algebras, but also to non-Gorenstein ones (see e.g. [B2], [BK],
[BR], [C], [CPST], [Rin], [T], [Y]).

\vskip5pt

\subsection{} Throughout $\mathcal A$ is an
abelian category with enough projective objects unless stated
otherwise, $\mathcal P(\mathcal A)$ is the full subcategory of
projective objects, and  $\mathcal{GP}(\mathcal A)$ is the full
subcategory of Gorenstein-projective objects. Replacing $\mathcal
P(\mathcal A)$ by a full additive subcategory $\mathcal C$, one
similarly define a category $\mathcal G(\mathcal C)$, and  $\mathcal
G^2(\mathcal C): = \mathcal G(\mathcal G(\mathcal C))$. S.
Sather-Wagstaff, T. Sharif and D. White proved that if $\mathcal{C}$
is self-orthogonal, then $\mathcal{G}^2(\mathcal{C}) =
\mathcal{G}(\mathcal{C})$; and they proposed the question: Whether
$\mathcal{G}^2(\mathcal{C}) = \mathcal{G}(\mathcal{C})$ holds for an
arbitrary $\mathcal{C}$? See [SWSW, 4.10, 5.8]. Recently  Z. Y.
Huang ([Huang, 4.1]) answered this question affirmatively. This
shows that $\mathcal{G}(\mathcal{C})$, in particular
$\mathcal{GP}(\mathcal A)$, has a strong stability. We give a slight
generalization of this stability by a different method. For details
please see 2.1 and Theorem \ref{stability}.

\vskip5pt

\subsection{} Throughout $\A$-mod is the category of finitely generated left modules of Artin
algebra $\A$, $\mathcal{P}(\A)$ the full subcategory of projective
$\A$-modules, and $\mathcal{GP}(\A)$ the full subcategory of
Gorenstein-projective $\A$-modules. Clearly $\mathcal{P}(\A)
\subseteq \mathcal{GP}(\A)$. If $\mathcal{GP}(\A) =
\mathcal{P}(\A)$, then $\A$ is called {\rm CM}-{\it free}. Algebras
of finite global dimension are CM-free ([EJ2, 10.2.3]), and there
are also many CM-free algebras of infinite global dimension ([C]).
If $\mathcal{GP}(\A)$ has only finitely many isoclasses of
indecomposable objects, then $\A$ is called {\rm CM}-{\it finite}.
In this case, let $E_1, \cdots, E_n$ be all the pairwise
non-isomorphic indecomposable Gorenstein-projective modules, and
${\rm \bf Aus}(\A)$ {\it the relative Auslander algebra} $\End_\A
(\bigoplus \limits_{1\le i\le n}E_i)^{op}$. CM-finiteness and
CM-freeness have been recently studied for examples in [B2], [B3],
[C], [CY], [GZ], [LZ], \cite{Rin}. In fact, CM finiteness
characterizes the simple hypersurface singularity ([BGS], [CPST],
[Kn]).

\vskip5pt

A CM-finite algebra $\A$ is Gorenstein if and only if ${\rm
gl.dim}{\rm \bf Aus}(\A) < \infty$ ([LZ]; [B3]). Thus in this case
${\rm \bf Aus}(\A)$ is CM-free. A CM-finite Gorenstein algebra $\A$
which is not CM-free has a non-zero singularity category, but the
singularity category of ${\rm \bf Aus}(\A)$ becomes zero. So, this
is a kind of categorical resolution of singularities ([Kuz]).
However, by the recent work of C. M. Ringel [Rin] there are many
CM-finite algebras which are non-Gorenstein and not CM-free (see
2.3). So it is natural to study ${\rm \bf Aus}(\A)$ also for such
algebras. We prove that for an arbitrary CM-finite algebra $\A$,
${\rm \bf Aus}(\A)$ is always CM-free (Theorem \ref{free}). It is a
Gorenstein version of Auslander's theorem ([ARS, p.215]). For those
CM-finite algebras which are non-Gorenstein and not CM-free, Theorem
\ref{free} means a categorical resolution of Gorenstein
singularities.

\vskip5pt

Put $\Omega$ to be the class of CM-finite algebras, and $\Theta$ the
class of CM-free algebras. Theorem \ref{free} implies that there is
a surjective map ${\rm \bf Aus}: \ \Omega \longrightarrow \Theta$;
moreover, it sends CM-finite Gorenstein algebras to algebras of
finite global dimension, and sends CM-finite non-Gorenstein algebras
to CM-free non-Gorenstein algebras. Graphically, we have

\begin{center}
\includegraphics[width = 150mm]{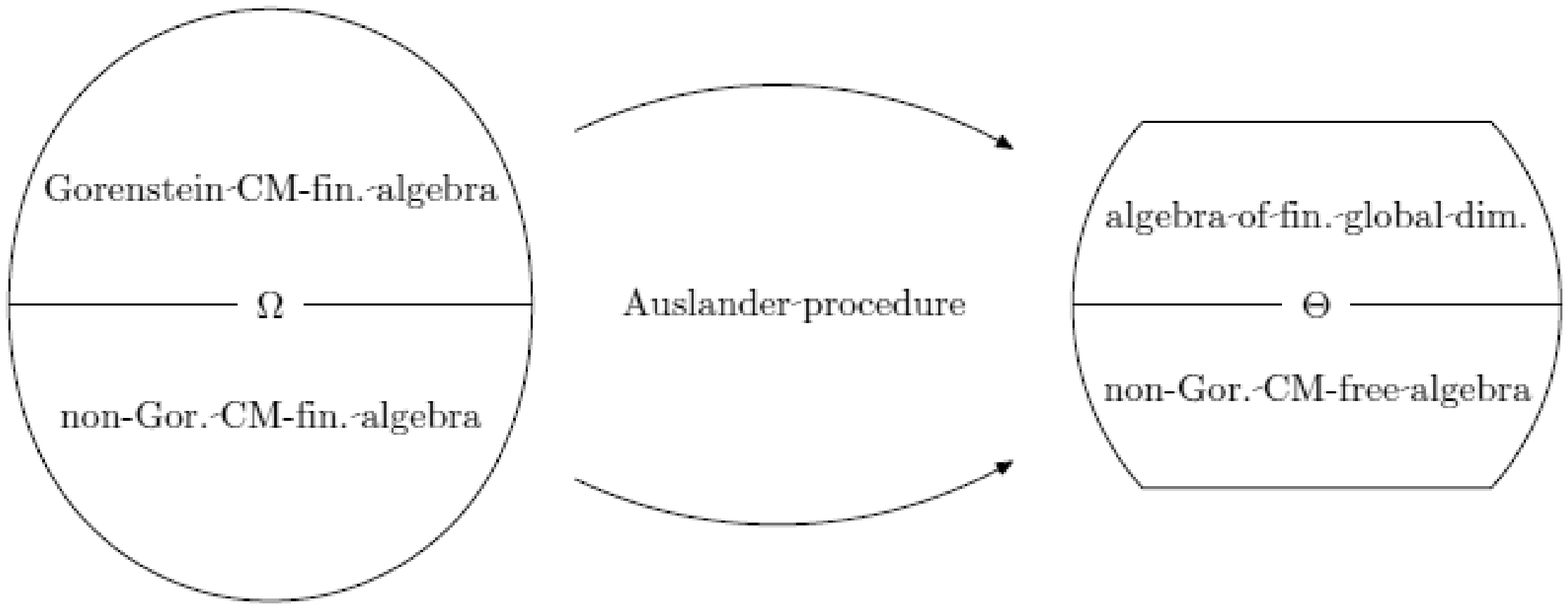}
\end{center}

\subsection{} An important feature is that $\mathcal {GP}(\mathcal A)$ is a
Frobenius category ([Ke]) with relative projective-injective objects
being projective objects ([B2]), and hence the stable category
$\underline{\mathcal {GP}(\mathcal A)}$ of $\mathcal {GP}(\mathcal
A)$ modulo $\mathcal {P}(\mathcal A)$ is triangulated (see D. Happel
[Hap1, p.16]).

\vskip5pt

The singularity category $D^b_{sg}(\mathcal A)$ of $\mathcal A$ is
defined as the Verdier quotient $D^b(\mathcal A)/K^{b}(\mathcal
P(\mathcal A))$ (R.-O. Buchweitz [Buch], D. Orlov [O]), where
$D^b(\mathcal A)$ is the bounded derived category of $\mathcal A$,
and $K^b(\mathcal P(\mathcal A))$ is the thick subcategory
consisting of bounded complexes of projective objects. It measures
how far $\mathcal A$ is from smoothness. Buchweitz's Theorem ([Buch,
4.4.1]; also [Hap2, 4.6]) says that there is a triangle-embedding
$F: \underline {\mathcal{GP}(\mathcal A)} \hookrightarrow
D^b_{sg}(\mathcal A)$, and that if each object of $\mathcal A$ is of
finite Gorenstein-projective dimension, then $F$ is an equivalence
(it is stated for $R$-mod, but it holds also for $\mathcal A$). The
converse is also true. See A. Beligiannis [B1, 6.9(8)], [BJO], and
 S. J. Zhu [Zhu].

\vskip5pt

Following P. A. Bergh, D. A. Jorgensen and S. Oppermann [BJO], the
Gorenstein defect category of $\mathcal A$ is defined as the Verdier
quotient $D^b_{{\rm defect} }(\mathcal A): = D^b_{sg}(\mathcal
A)/{\rm Im}F.$
 It measures how far $\mathcal A$ is from Gorensteinness since
$D^b_{\rm defect}(\mathcal A)=0$ if and only if each object of
$\mathcal A$ is of finite Gorenstein-projective dimension.

\vskip5pt

Another aim of this paper is to describe $D^b_{{\rm defect}
}(\mathcal A)$.

\vskip5pt

\subsection{} Note that $D^b(\mathcal A)$ can
be interpreted as $K^{-, b}(\mathcal P(\mathcal A))$. This makes
$D^b(\mathcal A)$ more accessible. If $\mathcal{A}$ is {\rm
CM}-contravariantly finite, we can describe $D^b(\mathcal A)$ also
via Gorenstein-projective objects, i.e., there is a
triangle-equivalence
$$D^b(\mathcal A)\cong {K^{-, gpb}(\mathcal{GP}(\mathcal A))}/{K^{b,
ac}(\mathcal{GP}(\mathcal A))},$$  where ${K^{-,
gpb}(\mathcal{GP}(\mathcal A))}$ is introduced in [GZ]. For details
please see 2.2 and Theorem \ref{DB}. It plays an important role in
our investigation on $D^b_{{\rm defect} }(\mathcal A)$.

\vskip5pt

\subsection{}  By introducing the category $K^{-, \ b}_\mathcal G(\mathcal P(\mathcal A))$, we describe $D^b_{{\rm defect} }(\mathcal
A)$ as
$$D^b_{{\rm defect}}(\mathcal A) \cong K^{-, \ b}(\mathcal P(\mathcal A))/K^{-, \ b}_\mathcal G(\mathcal P(\mathcal A));$$
and if $\mathcal A$ is {\rm CM}-contravariantly finite, then we also
have
$$D^b_{{\rm defect}}(\mathcal A) \cong K^{-, gpb}(\mathcal {GP}
(\mathcal A))/K^{b}(\mathcal {GP} (\mathcal A)).$$ For details
please see Theorem \ref{descriptionofdeft}. This implies that in
this condition $D^b_{{\rm defect}}(\mathcal A)$ is completely
controlled by $\mathcal {GP} (\mathcal A)$.  As an application, if
$\A$ is CM-finite, then we get a triangle-equivalence (Corollary
\ref{dif})
$$D^b_{{\rm defect}}(\A) \cong  D_{sg}^b( {\rm \bf Aus}(\A)).$$
Thus the Gorenstein defect category of a CM-finite algebra can be
reduced to the singularity category of a CM-free algebra.

\vskip5pt

\subsection{} Suppose that $\mathcal A$ is
{\rm CM}-contravariantly finite with $\mathcal{GP}(\mathcal A) \cong
\mathcal{GP}(\mathcal A')$.  One can not expect that $\mathcal A'$
is also {\rm CM}-contravariantly finite. However, $D^b(\mathcal A')$
and $D^b_{defect}(\mathcal A')$ share the same descriptions via
$\mathcal{GP}(\mathcal A')$. So $D^b_{defect}(\mathcal A) \cong
D^b_{defect}(\mathcal A').$ Please see Theorem \ref{frem}.

\vskip5pt

\subsection{} \ The paper is organized as follows. In
Section 2 we give necessary preliminaries and notations. In Section
3 we prove a version of the stability of $\mathcal G(\mathcal C)$.
Sections 4 is devoted to Theorem \ref{free}. Sections 5 and 6 are to
describe $D^b(\mathcal A)$ and $D^b_{{\rm defect} }(\mathcal A)$ via
Gorenstein-projective objects. Section 7 is to prove Theorem
\ref{frem}.

\vskip5pt

\section {\bf Preliminaries and notations}
\subsection{Category $\mathcal G(\mathcal C)$ and Gorenstein-projective objects} \ \ Let
$\mathcal{C}$ be a full additive subcategory of an abelian category
$\mathcal{A}$ (not necessarily has enough projective objects) which
is closed under isomorphisms. A complex $X^\bullet$ over
$\mathcal{A}$ is {\it $\mathcal{C}$-exact} (resp. {\it
$\mathcal{C}$-coexact}) if $\Hom_\mathcal A(C, X^\bullet)$ is exact
(resp. $\Hom_\mathcal A(X^\bullet, C)$ is exact) for each $C\in
\mathcal{C}$. A $\mathcal{C}$-exact and $\mathcal{C}$-coexact
complex is said to be {\it $\mathcal{C}$-biexact}.

Let $L$ be an object of $\mathcal{A}$. An exact complex $X^\bullet =
(X^i, d^i)$ over $\mathcal{A}$ is {\it a complete
$\mathcal{C}$-resolution of $L$}, if each $X^i \in \mathcal{C}$ and
$X^\bullet$ is $\mathcal{C}$-biexact, such that $L \cong \Ima d^0$.
Define $\mathcal{G}(\mathcal{C})$ to be the full subcategory of
$\mathcal{A}$ consisting of the objects which admit complete
$\mathcal{C}$-resolutions. Define $\mathcal{G}^2(\mathcal{C}): =
\mathcal{G}(\mathcal G(\mathcal{C}))$  ([SWSW, 4.1]). Clearly
$\mathcal C \subseteq \mathcal{G}(\mathcal{C}) \subseteq
\mathcal{G}^2(\mathcal{C}).$

Let $\mathcal A$ be  an abelian category with enough projective
objects. Taking $\mathcal C$ to be $\mathcal P(\mathcal A)$, we get
the notion of {\it a complete projective resolution}. Let
${\mathcal{GP}(\mathcal A)}$ denote the category
$\mathcal{G}(\mathcal P(\mathcal A))$, whose objects are called {\it
Gorenstein-projective objects} ([AB], [EJ1]).
\subsection{CM-contravariantly finite abelian categories}
Let $\mathcal{A}$ be an abelian category (not necessarily has enough
projective objects),  $\mathcal{C}$ a full subcategory of $\mathcal
A$, and $X\in\mathcal A$. A morphism  $f: C \longrightarrow X$ with
$C\in \mathcal{C}$ is {\it a right $\mathcal{C}$-approximation} of
$X$, if $\Hom_{\mathcal A}(C', f): \Hom_{\mathcal A}(C', C)
\longrightarrow \Hom_{\mathcal A}(C', X)$ is surjective for each $C'
\in \mathcal{C}$.  If each object $X\in \mathcal A$ admits a right
$\mathcal{C}$-approximation, then $\mathcal{C}$ is said to be {\it
contravariantly finite in $\mathcal A$} (\cite{AR}).

\vskip5pt

Let $\mathcal{A}$ be an abelian category with enough projective
objects. For short, we say that $\mathcal A$ is {\it {\rm
CM}-contravariantly finite},  if $\mathcal{GP}(\mathcal A)$ is
contravariantly finite in $\mathcal A$. If each object of
$\mathcal{A}$ has a finite Gorenstein-projective dimension, then
$\mathcal A$ is CM-contravariantly finite ([EJ2, 11.5.1], or H.Holm
[Hol, 2.10]. We stress that the proof in [Hol, 2.10] is stated for
module category over ring, but holds also for an abelian category
with enough projective objects).

An Artin algebra $\A$ is {\it {\rm CM}-contravariantly finite}, if
$\A$-mod is CM-contravariantly finite. Recall that $\A$ is {\it a
Gorenstein algebra}, if the injective dimension of $_\A\A$ is finite
and the injective dimension of $\A_\A$ is finite (in this case the
both are same, see [I]). Note that $\A$ is a Gorenstein algebra if
and only if each $\A$-module has finite Gorenstein-projective
dimension ([Hos]).  Thus a Gorenstein algebra is CM-contravariantly
finite. Also, clearly a CM-finite algebra is CM-contravariantly
finite. By [B2, Theorem 8.2${\rm (ix)}$] a virtually Gorenstein
algebra is CM-contravariantly finite. On the other hand, there
exists an Artin algebra which is not {\rm CM}-contravariantly finite
(see [BK], [Y]. Also [T]).

\subsection{CM-finite algebras which are non-Gorenstein and not CM-free}
There are many CM-finite algebras which are non-Gorenstein and not
CM-free. For examples, by [Rin, Proposition 5] the Nakayama algebras
with admissible sequences $(6,6,5)$, $(8, 8, 8, 7), \ (10, 10, 9,
10, 9)$, are such examples.

\subsection{Triangulated categories} \ Among the other
conditions, we emphasize that a triangulated subcategory $\mathcal
D$ of a triangulated category $\mathcal C$ is a full subcategory and
is closed under isomorphisms of $\mathcal C$. See A. Neeman [N]. For
a triangle functor $F: \mathcal A \longrightarrow \mathcal B$
between triangulated categories, let ${\rm Im}F$ denote the full
subcategory of $\mathcal B$ consisting of the objects which are
isomorphic to $F(X)$ with $X\in\mathcal A$. If $F$ is full, then
${\rm Im}F$ is a triangulated subcategory of $\mathcal B$.
\begin{lem} \label{verdierquotient} \ {\rm (J. L. Verdier [V], \mbox{Corollary }  4\mbox{-}3)} \ Suppose $\mathcal{D}_1$ and $\mathcal{D}_2$
are triangulated subcategories of triangulated category
$\mathcal{C}$, and $\mathcal{D}_1$ is a subcategory of
$\mathcal{D}_2$. Then $\mathcal{D}_2/\mathcal{D}_1$ is a
triangulated subcategory of $\mathcal{C}/\mathcal{D}_1$, and there
is a triangle-equivalence
$(\mathcal{C}/\mathcal{D}_1)/(\mathcal{D}_2/\mathcal{D}_1)\cong
\mathcal{C}/\mathcal{D}_2$.
\end{lem}
\begin{lem} \label{ff} \ {\rm ( J. Rickard \cite{Ric}, p. 446, line 1)}  A full triangle functor which sends
non-zero objects to non-zero objects is faithful.
\end{lem}
\subsection{Notations and convention} \ For convenience, we
list some categories mainly used in Section 5 - Section 7.

\vskip5pt

$K^{b}(\mathcal A) =$ the homotopy category of bounded complexes
over $\mathcal A$.

$K^{b, ac}(\mathcal A)=$ the homotopy category of bounded exact
complexes over $\mathcal A$.

$K^{-}(\mathcal A) =$ the homotopy category of upper bounded
complexes over $\mathcal A$.

$K^{-, ac}(\mathcal A)=$ the homotopy category of upper bounded
exact complexes over $\mathcal A$.

\vskip5pt

$K^{b}(\mathcal{P}(\mathcal A))=$ the homotopy category of bounded
complexes over $\mathcal P(\mathcal A)$.

$K^{-}(\mathcal{P}(\mathcal A))=$ the homotopy category of upper
bounded complexes over $\mathcal P(\mathcal A)$.

$K^{-, b}(\mathcal{P}(\mathcal A))=$ the homotopy category of upper
bounded complexes over $\mathcal P(\mathcal A)$, with only finitely
many non-zero cohomologies.

\vskip5pt

$K^{b}(\mathcal{GP}(\mathcal A))=$ the homotopy category of bounded
complexes over $\mathcal{GP}(\mathcal A)$.

$K^{b, ac}(\mathcal{GP}(\mathcal A))=$ the homotopy category of
bounded exact complexes over $\mathcal{GP}(\mathcal A)$. Here
``exact" means exact as a complex over $\mathcal A.$

$K^{-}(\mathcal{GP}(\mathcal A))=$ the homotopy category of upper
bounded complexes over $\mathcal{GP}(\mathcal A)$.

\vskip5pt

$D^{b}(\mathcal A)=$ the derived category of bounded complexes over
$\mathcal A$,  i.e., the Verdier quotient $K^{b}(\mathcal A)/K^{b,
ac}(\mathcal A)$.

$\langle \mathcal{GP}(\mathcal A)\rangle =$ the triangulated
subcategory of $D^b(\mathcal A)$ generated by $\mathcal{GP}(\mathcal
A)$, i.e., the smallest triangulated subcategory of $D^b(\mathcal
A)$ containing $\mathcal{GP}(\mathcal A)$.

$D^{-}(\mathcal A)=$ the derived category of upper bounded complexes
over $\mathcal A$, i.e., the Verdier quotient $K^{-}(\mathcal
A)/K^{-, ac}(\mathcal A)$.

$\underline{\mathcal {GP}(\mathcal A)}=$ the stable category of
$\mathcal {GP}(\mathcal A)$ modulo $\mathcal {P}(\mathcal A)$.

\vskip5pt

{\bf Convention.} \ Let $\mathcal C$ be a triangulated category and
$\mathcal D$  a full additive subcategory. We say that $\mathcal D$
is a triangulated subcategory of $\mathcal C$, if the isomorphism
closure $\widetilde{\mathcal D}$ of $\mathcal D$ is a triangulated
subcategory of $\mathcal C$. In this case,  we do not distinguish
between $\mathcal D$ and $\widetilde{\mathcal D}$ (if no substantial
difficulties occur): \ for example, the Verdier quotient $\mathcal
C/\mathcal D$ always means $\mathcal C/\widetilde{\mathcal D}$. This
convention occurs in many places in this paper.

\vskip5pt

\subsection{Singularity categories} The Verdier quotient $D^b(\mathcal A)/K^{b}(\mathcal
P(\mathcal A))$ is called {\it the singularity category of $\mathcal
A$} and denoted by $D^b_{sg}(\mathcal A)$ (see [O]; or {\it the
stabilized  derived category} in [Buch]). Then  $D^b_{sg}(\mathcal
A)=0$ if and only if each object of $\mathcal A$ has finite
projective dimension. By the canonical triangle-equivalence $\rho:
K^{-, b}(\mathcal P(\mathcal A)) \longrightarrow D^b(\mathcal A)$,
there is a triangle-equivalence $D^b_{sg}(\mathcal A)\cong K^{-,
b}(\mathcal P(\mathcal A))/ K^{b}(\mathcal P(\mathcal A)).$

\vskip5pt

\subsection{Buchweitz Theorem and Gorenstein defect categories}
Consider the composition of the embedding $\mathcal {GP}(\mathcal A)
\hookrightarrow  D^b(\mathcal A)$ and the localization functor
$D^b(\mathcal A) \longrightarrow D^b(\mathcal A)/K^b(\mathcal
P(\mathcal A)).$ It induces a functor
$$F:  \underline {\mathcal{GP}(\mathcal A)} \longrightarrow D^b_{sg}(\mathcal
A) : = D^b(\mathcal A)/K^b(\mathcal P(\mathcal A))\eqno(2.1)$$ which
sends a Gorenstein-projective object $G$ to the stalk complex of $G$
at degree $0$. By the triangle-equivalence $\overline{\rho}: K^{-,
b}(\mathcal P(\mathcal A))/ K^{b}(\mathcal P(\mathcal A))\cong
D^b(\mathcal A)/K^b(\mathcal P(\mathcal A))$ one gets a functor
$$\overline{F}:  \underline {\mathcal{GP}(\mathcal A)} \longrightarrow K^{-, b}(\mathcal P(\mathcal A))/ K^{b}(\mathcal P(\mathcal A))\eqno(2.2)$$
such that the diagram
\[\xymatrix {\underline {\mathcal{GP}(\mathcal A)} \ar@{=}[d]\ar[r]^-{\overline F} & K^{-, b}(\mathcal P(\mathcal A))/ K^{b}(\mathcal P(\mathcal A)) \ar[d]^-{\overline{\rho}}
\\ \underline {\mathcal{GP}(\mathcal A)} \ar[r]^-{F} & D^b(\mathcal A)/K^b(\mathcal P(\mathcal
A))}\eqno(2.3)\] commutes. For a Gorenstein-projective object $G$,
let $(P^\bullet, d)$ be a complete projective resolution such that
${\rm Im} d^{0}\cong G$. Then
$$\overline{F}(G) = P^\bullet_{\le 0}\in K^{-, b}(\mathcal
P(\mathcal A))/ K^{b}(\mathcal P(\mathcal A)),$$ where
$P^\bullet_{\le 0}$ is the brutal truncation $\cdots \longrightarrow
P^{-1} \stackrel {d^{-1}}\longrightarrow P^0 \longrightarrow 0$ of
$P^\bullet$.

\vskip5pt

Following [BJO], the Gorenstein defect category of $\mathcal A$ is
defined as the Verdier quotient
$$D^b_{{\rm defect} }(\mathcal A): = D^b_{sg}(\mathcal A)/{\rm Im}F = (D^b(\mathcal A)/K^b(\mathcal P(\mathcal
A)))/{\rm Im}F.$$ For an Artin algebra $\A$, the Gorenstein defect
category $D^b_{\rm defect}(\A)$ of $\A$ is defined to be $D^b_{\rm
defect}(\A\mbox{-}{\rm mod})$. Then $D^b_{\rm defect}(\A) = 0$ if
and only if $\A$ is Gorenstein. By $(2.3)$ the diagram
\[\xymatrix {{\rm Im}\overline{F}\ar[r]\ar[d] & K^{-, b}(\mathcal P(\mathcal A))/ K^{b}(\mathcal P(\mathcal A)) \ar[d]^-{\overline{\rho}}
\\ {\rm Im}F\ar[r] & D^b(\mathcal A)/K^b(\mathcal P(\mathcal
A))}\eqno(2.4)\] commutes, where the horizontal functors are
embeddings, and the vertical functors are triangle-equivalences. It
follows that there is a triangle-equivalence
$$D^b_{{\rm defect} }(\mathcal A)\cong (K^{-, b}(\mathcal P(\mathcal A))/ K^{b}(\mathcal P(\mathcal A)))/{\rm Im}\overline{F}.\eqno(2.5)$$

\vskip5pt

\section {\bf Stability of Gorenstein categories}

Throughout $\mathcal{A}$ is an abelian category (not necessarily has
enough projective objects), and  $\mathcal{C}$ is a full additive
subcategory of $\mathcal{A}$, which is closed under isomorphisms. By
[SWSW, 4.10] and [Huang, 4.1] $\mathcal{G}^2(\mathcal{C}) =
\mathcal{G}(\mathcal{C})$. This shows that the category
$\mathcal{G}(\mathcal{C})$ has a strong stability. By using a
different method, we have a little generalization of this stability.

\vskip5pt

\begin{thm} \label{stability} \ Let $\mathcal{A}$ be an abelian category, and
$\mathcal{C}$ a full additive subcategory of $\mathcal{A}$ which is
closed under isomorphisms. Let  $X^\bullet = (X^i, d^i)$ be an exact
complex with $X^i  \in \mathcal{G}(\mathcal{C})$ for all $i\in\Bbb
Z$. If $X^\bullet$ is $\mathcal{C}$-biexact, then $\Ima d^i \in
\mathcal{G}(\mathcal{C})$ for all $i\in\Bbb Z$.
\end{thm}

If $X^\bullet$ in Theorem \ref{stability} is required to be
$\mathcal G({\mathcal C})$-biexact, then it is exactly
$\mathcal{G}^2(\mathcal{C}) = \mathcal{G}(\mathcal{C})$.

\subsection{} We need the following fact.

\begin{lem} \label{push-out} \ Suppose there is a commutative diagram in $\mathcal{A}$ with exact rows $\delta$ and $\eta$:
$$ \xymatrix{
\delta : &0  \ar[r] &X \ar[r]^f \ar[d]^{\alpha} &Y \ar[d]^{\beta} \ar[r]^g &Z \ar[r] \ar@{=}[d] &0 \\
\eta :   &0   \ar[r] &X^\prime \ar[r]^{f^\prime}         &Y^\prime
\ar[r]^{g^\prime}                   &Z   \ar[r]              &0.  }
$$ Denote by $\Delta$ the corresponding short exact sequence $$ 0
\longrightarrow X \xlongrightarrow{\left(  \begin{smallmatrix} - f \\
\alpha \end{smallmatrix}\right)} Y \oplus X^\prime
\xlongrightarrow{\left(\begin{smallmatrix} \beta, & f^\prime
\end{smallmatrix}\right)} Y^\prime \longrightarrow 0.$$ Then we have

\vskip5pt

$(i)$ \ $\delta$ is $\mathcal C$-exact if and only if both $\eta$
and $\Delta$ are $\mathcal C$-exact.

\vskip5pt

$(ii)$  \ $\delta$ is $\mathcal C$-coexact if and only if  both
$\eta$ and $\Delta$ are $\mathcal C$-coexact.
\end{lem}

\noindent {\bf Proof.} $(i)$ \ We need to prove that  ${\rm
Coker}(\Hom_\mathcal A(C,  g)) = 0$ if and only if ${\rm Coker}
(\Hom_\mathcal A (C, g'))= 0$ and ${\rm Coker}(\Hom_\mathcal A(C,
(\beta, f^\prime))=0$, for each $C\in\mathcal C$. This can be seen
from the diagram chasing. However, for simplicity we use an argument
from triangulated category.

\vskip5pt

Regard $\delta$, $\eta$ and $\Delta$ as complexes in
$K^b(\mathcal{A})$, and the above commutative diagram as a morphism
$h$ from $\delta$ to $\eta$. Then the mapping cone ${\rm Cone}(h)$
is the complex $$ 0
\longrightarrow X \xlongrightarrow{\left(  \begin{smallmatrix} - f \\
\alpha \end{smallmatrix}\right)} Y \oplus X^\prime
\xlongrightarrow{\left(\begin{smallmatrix} -g & 0\\ \beta & f^\prime
\end{smallmatrix}\right)} Z\oplus Y^\prime \xlongrightarrow {(1, \  g')}Z \longrightarrow 0.$$
Then there is a homotopy equivalence $\Delta[1] \cong {\rm
Cone}(h)$, see the commutative diagram

\[\xymatrix{
0
\ar[r] &X \ar@{=}[d]  \ar[r]^{\left(  \begin{smallmatrix} - f \\
\alpha \end{smallmatrix}\right)} &Y \oplus X^\prime \ar[d]_{-{\rm
Id}} \ar[r]^{\left(\begin{smallmatrix} -g & 0\\ \beta & f^\prime
\end{smallmatrix}\right)} &Z\oplus Y^\prime  \ar[d]^{(0,\ 1)}   \ar[r]^-{(1, \  g')} &Z
\ar[d] \ar[r]
&0\\
0 \ar[r] & X \ar[r]^{\left(  \begin{smallmatrix} f \\
-\alpha \end{smallmatrix}\right)}  \ar@{=}[d]  & Y \oplus X^\prime
\ar[r]^{\left(\begin{smallmatrix} -\beta, & -f^\prime
\end{smallmatrix}\right)} \ar[d]_{-{\rm Id}}  & Y^\prime \ar[r] \ar[d]^-{\left(  \begin{smallmatrix} -g' \\
1 \end{smallmatrix}\right)}  & 0 \ar[r] \ar[d] &
0\\
0
\ar[r] &X \ar[r]^-{\left(  \begin{smallmatrix} - f \\
\alpha \end{smallmatrix}\right)} &Y \oplus X^\prime
\ar[r]^-{\left(\begin{smallmatrix} -g & 0\\ \beta & f^\prime
\end{smallmatrix}\right)} &Z\oplus Y^\prime \ar[r]^-{(1, \  g')} &Z \ar[r]
&0. }\]

\vskip5pt \noindent Thus we have a distinguished triangle $\delta
\longrightarrow \eta \longrightarrow \Delta[1] \longrightarrow
\delta[1]$ in $K^b(\mathcal{A})$. Applying the cohomological functor
$\Hom_{K^b(\mathcal A)}(C, -)$ to this distinguished triangle we get
the following long exact sequence of abelian groups
\begin{align*}\cdots \longrightarrow \Hom_{K^b(\mathcal A)}(C, \eta[1])
& \longrightarrow \Hom_{K^b(\mathcal A)}(C, \Delta[2])
\longrightarrow \Hom_{K^b(\mathcal A)}(C, \delta[2])\\&
\longrightarrow \Hom_{K^b(\mathcal A)}(C, \eta[2])\longrightarrow
\Hom_{K^b(\mathcal A)}(C, \Delta[3])\longrightarrow
\cdots.\end{align*}  Using the formula $\Hom_{K^b(\mathcal A)}(C,
\eta[n]) = {\rm H}^n\Hom_{\mathcal A}(C, \eta),$ the above exact
sequence read is
\begin{align*}\cdots \longrightarrow {\rm H}^1\Hom_{\mathcal
A}(C, \eta) & \longrightarrow {\rm H}^2\Hom_{\mathcal A}(C, \Delta)
\longrightarrow {\rm H}^2\Hom_{\mathcal A}(C, \delta)\\&
\longrightarrow {\rm H}^2\Hom_{\mathcal A}(C, \eta)\longrightarrow
{\rm H}^3\Hom_{\mathcal A}(C, \Delta)\longrightarrow
\cdots.\end{align*} That is we have the exact sequence
$$0
\longrightarrow {\rm Coker}(\Hom_{\mathcal A}(C, (\beta, f^\prime))
\longrightarrow {\rm Coker}(\Hom_{\mathcal A}(C, g)) \longrightarrow
{\rm Coker} (\Hom_{\mathcal A}(C, g')) \longrightarrow 0.$$ This
prove $(i)$.

\vskip5pt

$(ii)$ can be similarly proved. $\s$

\vskip5pt

\subsection{} We also need the following technical lemma.

\vskip5pt

\begin{lem} \label{induction}

Let $\delta: 0 \longrightarrow X_1 \xlongrightarrow{f} U
\xlongrightarrow{g} X_2 \longrightarrow 0 $ be an exact sequence
with $U \in \mathcal{G}(\mathcal{C})$ such that $\delta$ is
$\mathcal{C}$-biexact. Let $\eta: 0 \longrightarrow X_1
\xlongrightarrow{u} Y \xlongrightarrow{v} V \longrightarrow 0$ be an
exact sequence with $V \in \mathcal{G}(\mathcal{C})$ such that
$\eta$ is $\mathcal{C}$-coexact. Then we have a commutative diagram

\vskip5pt

\[ \xymatrix{
\delta: &0  \ar[r] &X_1 \ar[r]^f \ar@{_{(}->}[d]^{u} &U \ar[d] \ar[r]^g &X_2 \ar[r] \ar@{_{(}->}[d]^{u^\prime} &0 \\
\delta^\prime: &0   \ar[r] &Y \ar[r]         &C \ar[r]
&Z   \ar[r] &0 } \] such that

\vskip5pt

$(i)$ \ $\delta^\prime$ is  exact  with $C \in \mathcal{{C}}$,  and
$\delta^\prime$ is $\mathcal{C}$-biexact;

\vskip5pt

$(ii)$ \ $u^\prime$ is a monomorphism  with   $\Cok{u^\prime} \in
\mathcal{G}(\mathcal{{C}})$, and $\Hom_{\mathcal A}(u^\prime, \mathcal{C})$ is
surjective.
\end{lem}

\noindent {\bf Proof.} \ Step 1. \  Since by assumption $U, V \in
\mathcal{G}(\mathcal{C})$, by definition there exists two exact
sequences $\epsilon_1: 0 \longrightarrow U \xlongrightarrow{a} C_1
\xlongrightarrow{b} L_1 \longrightarrow 0$ and $\epsilon_2: 0
\longrightarrow V \xlongrightarrow{c} C_2 \xlongrightarrow{d} L_2
\longrightarrow 0$ with $C_1, C_2\in \mathcal{C}, L_1, L_2 \in
\mathcal{G}(\mathcal{C})$,  such that both $\epsilon_1$ and
$\epsilon_2$ are $\mathcal{C}$-biexact. Since by assumption
$\Hom_{\mathcal A}(\eta, C_1)$ is exact,  by considering $af\in\Hom_{\mathcal A}(X_1, C_1)$ we
see that there exists morphisms $e$ and $e^\prime$ such that the
following  diagram commutes:
$$ \xymatrix{ \eta: & 0 \ar[r] &X_1 \ar[d]^f  \ar[r]^u &Y   \ar@{-->}[d]^e  \ar[r]^v &V  \ar@{-->}[d]^{e^\prime} \ar[r] &0\\ &  0
\ar[r] &U \ar[r]^a &C_1  \ar[r]^b  &L_1 \ar[r] &0.}$$
\vskip5pt \noindent Put $\alpha =  \left(\begin{smallmatrix} cv \\
e
\end{smallmatrix} \right): Y\longrightarrow C_2\oplus C_1, \ i =
\left(\begin{smallmatrix} 0 \\ a \end{smallmatrix} \right):
U\longrightarrow C_2\oplus C_1$,
$\pi = \bigl(\begin{smallmatrix} -1 & 0 \\
0 & b \end{smallmatrix} \bigr): C_2\oplus C_1\longrightarrow
C_2\oplus L_1,  \ x =  \left(\begin{smallmatrix} -c \\  e'
\end{smallmatrix} \right): V\longrightarrow C_2\oplus L_1$. Then
there exist morphisms $\beta, u', v', y$ \ such that the following
diagram commutes:
$$ \xymatrix{ &
& & \eta      &\Xi    &\eta^\prime \\
& & & 0 \ar[d] &  0 \ar[d]  &0 \ar@{-->}[d]\\
& \delta: &0  \ar[r] &X_1 \ar[d]^u \ar[r]^-f &U \ar[d]^{i} \ar[r]^-g &X_2 \ar@{-->}[d]^{u^\prime} \ar[r] &0\\
& \delta^\prime: &0 \ar[r] &Y \ar[d]^{v}  \ar[r]^-{\alpha}
& C_2 \oplus C_1 \ar[r]^-\beta \ar[d]^{\pi}  & Z  \ar[r]  \ar@{-->}[d]^{v^\prime}   &0 & & (*)\\
& \Delta: &0 \ar[r] & V   \ar[r]^-{x} \ar[d]  & C_2 \oplus L_1 \ar[r]^-y \ar[d]  & L  \ar[r] \ar[d] & 0 \\
&&& 0  &  0 &0} $$ where $Z = \Cok \alpha$ and $L = \Cok x$. Since
$c$ is a monomorphism, so is $x$. Also, the middle column $\Xi$ \ is
exact since $\epsilon_1$ is exact. Applying Snake Lemma to the left
two columns we know that $\alpha$ is a monomorphism and that the
right column $\eta^\prime$ is exact. In particular, $u^\prime$ is a
monomorphism.

We will prove  that the upper two rows  of $(*)$ are what we need.

\vskip5pt

Step 2. Write $ y = (l, m)$. Observe that $\Delta$ is exact means
that there is the following commutative diagram with exact rows
(note that $L$ is the push-out of $e'$ and $c$):
$$ \xymatrix{
\epsilon_2:  &0  \ar[r] &V \ar[r]^c \ar[d]^{e^\prime} &C_2\ar[d]^l \ar[r]^d &L_2 \ar[r] \ar@{=}[d] &0 \\
\zeta:    &0   \ar[r] &L_1 \ar[r]^m         &L \ar[r]^n &L_2 \ar[r]
&0  } $$

\noindent Since  $\epsilon_2$ is $\mathcal{C}$-biexact, it follows
from Lemma \ref{push-out} that both $\zeta$ and $\Delta$ are
$\mathcal{C}$-biexact.

Since $L_1, L_2 \in \mathcal{G}(\mathcal{C})$, and $\zeta$ is
$\mathcal{C}$-biexact, it follows from Proposition 4.4 of
\cite{SWSW} that $L \in \mathcal{G}(\mathcal{C})$, thus $\Cok u' = L
\in \mathcal{G}(\mathcal{C})$.

\vskip5pt

Step 3.  $\forall \ C \in \mathcal{C}$, applying $\Hom_{\mathcal A}(C, -)$  to the right two columns of $(\ast)$, we get the following commutative diagram
with exact rows
$$ \xymatrix{
& &\Hom_{\mathcal A}(C, U) \ar[r] \ar[d]^{\Hom_{\mathcal A}(C, g)}
&\Hom_{\mathcal A}(C, C_2\oplus C_1)\ar[d]^{\Hom_{\mathcal A}(C,
\beta)} \ar[r]^{(C, \pi)} &\Hom_{\mathcal A}(C, C_2\oplus C_1)
\ar[r] \ar[d]^{\Hom_{\mathcal A}(C, y)} &0 \\
&0   \ar[r] &\Hom_{\mathcal A}(C, X_2)\ar[r] &\Hom_{\mathcal A}(C,
Z) \ar[r] &\Hom_{\mathcal A}(C, L).} $$ (Note that $\Hom_{\mathcal
A}(C, \pi)$ is surjective since  $\Hom_{\mathcal A}(C, \epsilon_1)$
is exact.) Also $\Hom_{\mathcal A}(C, y)$ and $\Hom_{\mathcal A}(C,
g)$  are surjective, by Snake Lemma $\Hom_{\mathcal A}(C, \beta)$ is
surjective. Thus $\delta^\prime$ is $\mathcal{C}$-exact.

\vskip5pt

Step 4.  $\forall \ C \in \mathcal{C}$, applying $\Hom_{\mathcal A}(-, C)$ to $\left(  \ast \right)$,
we have the following commutative diagram with exact rows and columns:

 $$ \xymatrix{
                     &         &(\eta^\prime, C)                                &(\Xi, C)                                            &(\eta, C)                 \\
                 &               &0      \ar[d]                             &0    \ar[d]                                   &0
 \ar[d]  \\
(\Delta, C):          &0    \ar[r]     &(L, C)      \ar[r]^-{(y, C)}
\ar[d]^{(v^\prime, C)}       &(C_2 \oplus L_1, C)
\ar[d]^{(\pi, C)}       \ar[r]^-{(x, C)}           &(V,C)    \ar[d]^{(v, C)}       \ar[r] &0   \\
(\delta^\prime, C):           &0    \ar[r]     &(Z, C)
\ar[d]^{(u^\prime, C)}    \ar[r]^-{(\beta, C)}
 &(C_2 \oplus C_1, C)       \ar[d]^{(i, C)}       \ar[r]^-{(\alpha, C)}            &(Y, C)   \ar[d]^{(u, C)}       \\
(\delta, C):    &0    \ar[r]     &(X_2, C)       \ar[r]^-{(g, C)}
&(U,  C)           \ar[r]^-{(f, C)}  \ar[d]
&(X_1,  C)    \ar[d]   \ar[r] &0 \\
              &          &       &0   &0
}$$
(Note that the exactness of $\Hom_{\mathcal A}(\delta, C)$ and $\Hom_{\mathcal A}(\eta, C)$ follows from assumptions; the exactness of $\Hom_{\mathcal A}(\Xi, C)$ follows from
the exactness of $\Hom_{\mathcal A}(\epsilon_1, C)$; and the exactness of $\Hom_{\mathcal A}(\Delta, C)$ follows from Step 2).
Applying Snake Lemma to the right two columns we see that both $\Hom_{\mathcal A}(u^\prime, C)$
and $\Hom_{\mathcal A}(\alpha, C)$ are surjective. This completes the proof.
$\s$

\vskip5pt

\subsection{Proof of Theorem \ref{stability}.} \ Without loss of generality, we only  prove
$\Ima d^0  \in \mathcal{G}(\mathcal{C})$. For each $i$, decompose
$d^i$ as $\xymatrix{X^i \ar@{->>}[r]^-{a^i}  &\Ima d^i
\ar@{^{(}->}[r]^{b^i} &X^{i+1}}$. We claim that there exists  an
exact sequence $C^+:  0 \longrightarrow \Ima d^0
\xlongrightarrow{v^0} C^{1}\xlongrightarrow{\gamma^{1}} C^{2}
\xlongrightarrow{\gamma^{2}} \cdots$ with $ C^i \in \mathcal{C}$ for
all $i \geq 1$, such that $C^+$ is $\mathcal{C}$-biexact.

In fact, applying Lemma \ref{induction} to the exact sequence
$0 \longrightarrow \Ima d^0  \xlongrightarrow {b^0} X^{1} \xlongrightarrow{a^{1}} \Ima d^{1}
\longrightarrow 0$ and the exact sequence  $0 \longrightarrow \Ima d^0 \xlongrightarrow{=} \Ima d^0\longrightarrow 0\longrightarrow 0$,
we get the following commutative diagram with exact rows
  \[\xymatrix{
 & 0 \ar[r] &\Ima d^0 \ar[r]^{b^{0}} \ar@{=}[d]&X^{1} \ar[d] \ar[r]^{a^{1}} & \Ima d^{1} \ar[r] \ar@{_{(}->}[d]^{w^{1}} &0\\
\xi_1: \ & 0 \ar[r] &\Ima d^0 \ar[r]^{v^0}                     &C^{1}           \ar[r]^{u^{1}}          &L^{1} \ar[r]     &0
 }\]
such that $C^{1}  \in \mathcal{C}$,  $w^1$ is a monomorphism with
$\Cok w^{1} \in \mathcal{G}(\mathcal{C})$, and that $\xi_1$ is
$\mathcal{C}$-biexact, and $\Hom_{\mathcal A}(w^{1}, \mathcal{C})$
is surjective. Now applying Lemma \ref{induction} to the exact
sequence $0 \longrightarrow \Ima d^{1} \xlongrightarrow{b^{1}} X^{2}
\xlongrightarrow{a^{2}} \Ima d^{2} \longrightarrow 0$ and the exact
sequence  $0 \longrightarrow \Ima d^1 \xlongrightarrow{w^1}
L^1\longrightarrow \Cok w^1\longrightarrow 0$, we get an exact
sequence $\xi_2: 0\longrightarrow L_1 \xlongrightarrow{v^1} C_2
\xlongrightarrow {u^2} L_3\longrightarrow 0$ and a monomorphism
$w^2: \Ima d^2 \hookrightarrow L^2$. Continuing this process and
putting $\xi_1, \xi_2, \cdots$ together we finally obtain an exact
sequence $C^+$ with $\gamma^i = v^iu^i, \ i\ge 1$, such that $C^+$
is $\mathcal{C}$-biexact (since each $\xi_i$ is
$\mathcal{C}$-biexact).

By duality there exists an exact sequence $C^-: \cdots
\longrightarrow C^{-1} \longrightarrow C^0 \longrightarrow \Ima{d^0}
\longrightarrow 0$ with $ C^i \in \mathcal{C}$ for $i \le 0$ such
that $C^-$ is $\mathcal{C}$-biexact (this can be done by applying
the above claim  to the opposite category $\mathcal{A}^{\rm op}$).
Now putting $C^-$ and $C^+$ together, we obtain a complete
$\mathcal{C}$-resolution of $\Ima{d^0}$, and this proves  $\Ima{d^0}
\in \mathcal{G}(\mathcal{C})$.  $\s$

\vskip5pt

\subsection{} Applying Theorem  \ref{stability} to $R$-mod, the category of
finitely generated left $R$-modules, we get

\begin{cor}\label{Rmod} Let $R$ be a left noetherian ring. Let $X^\bullet = (X^i, d^i)$ be an exact complex with $X^i
\in \mathcal{GP}(R)$ for all $i$. If \ $\Hom_{\A}(X^\bullet, R)$ is
still exact, then $\Ima d^i \in \mathcal{GP}(R)$ for all $i$.
\end{cor}

\section {\bf Relative Auslander algebras of CM-finite algebras}

Throughout this section $\A$ is an Artin algebra. All modules are
finitely generated.

\vskip5pt

\subsection{} Let $E$ be a Gorenstein-projective generator, i.e., $E$ is a
Gorenstein-projective $\A$-module such that $\A \in {\rm add} E$,
where ${\rm add} E$ is the smallest full additive subcategory of
$\A$-mod containing $E$ and closed under direct summands. For short,
we say that a complex $X^\bullet$ is {\it $E$-biexact}, if it is
$\add E$-biexact. Denote $\mathcal{G}(\add E)$ by $\mathcal{G}(E)$,
i.e.,
\begin{align*} \mathcal{G}(E) = \{X\in \A\mbox{-}{mod} \
 \mid \ & \exists \ \mbox{an exact complex} \ X^\bullet = (X^i, d^i) \
\text{with all} \  X^i \in \add E,  \\ & \ \mbox{such that} \
X^\bullet \ \mbox{is} \ E\mbox{-biexact}, \ \mbox{and} \  X\cong
\Ima d^0\}.\end{align*} By Corollary \ref{Rmod}$(ii)$ we have

\vskip5pt

\begin{cor} \label {GPE} \ Let $E$ be a
Gorenstein-projective generator. Then $\mathcal{G}(E) \subseteq
\mathcal{GP}(\A)$.
  \end{cor}
\subsection{} Recall the Yoneda philosophy. Given a
$\A$-module $E$, the  functor $\Hom_{\A}(E, -): \A\mbox{-}{\rm mod}
\longrightarrow \Gamma\mbox{-}{\rm mod}$ induces an equivalence
between $\add E$ and $\mathcal{P}(\Gamma)$, where $\Gamma = (\End_\A
E)^{op}$. Also, for each $E'\in\add E$ and each $X\in \A$-mod there
is an isomorphism
$$\Hom_\A (E', X) \longrightarrow \Hom_\Gamma (\Hom_{\A}(E, E'), \Hom_{\A}(E, X))$$
given by $f \mapsto \Hom_\A (E, f), \ \forall \ f\in \Hom_\A (E',
X)$  (cf. \cite{ARS}, p.33).

\vskip5pt

Now,  if $E$ is a generator of $\A$-mod (i.e., $\A\in\add E$), then
we can say more.

\vskip5pt

\begin{lem} \label{yoneda} \ Let $E$ be a generator of $\A$-mod, and $\Gamma = (\End E)^{op}$.
Then the functor $\Hom_{\A}(E, -): \A\mbox{-}{\rm mod}
\longrightarrow  \Gamma\mbox{-}{\rm mod}$ is fully faithful.
\end{lem}

\noindent {\bf Proof.}  Since $E$ is a generator, for any $\A$-module $X$ there is a surjective map
$E^m \twoheadrightarrow X$ for some positive integer $m$. This implies that
$\Hom_{\A}(E, -)$ is faithful. Let $X, Y \in \A$-mod, $f: \Hom_\A(E,
X) \longrightarrow  \Hom_\A(E, Y)$ be  a
  $\Gamma$-map. By taking right $\add E$-approximations,
we get exact sequences $T_1 \xlongrightarrow{u} T_0
\xlongrightarrow{\pi} X \longrightarrow 0$ and $T_1^\prime
\xlongrightarrow{u^\prime} T_0^\prime \xlongrightarrow{\pi^\prime} Y
\longrightarrow 0$ with $T_0, T_1, T_0', T_1'\in\add E$ (note that
since $E$ is a generator, it follows that $\pi$ and $\pi'$ are
surjective). Applying $\Hom_\A(E, -)$ we have the following diagram
with exact rows

  \[ \xymatrix{
  \Hom_{\A}(E, T_1) \ar[r]^{\Hom_{\A}(E,u)}  \ar@{-->}[d]^{f_1} &\Hom_{\A}(E, T_0) \ar[r]^{\Hom_{\A}(E,\pi)} \ar@{-->}[d]^{f_0} &\Hom_{\A}(E, X)
  \ar[r] \ar[d]^f &0\\
\Hom_{\A}(E, T_1^\prime) \ar[r]^{\Hom_{\A}(E,u^\prime)} &\Hom_{\A}(E, T_0^\prime)
\ar[r]^{\Hom_{\A}(E,\pi^\prime)}  &\Hom_{\A}(E, Y) \ar[r] &0.\\
  }\]
Since the two rows are respectively parts of projective resolutions
of $\Hom_{\A}(E, X)$ and $\Hom_{\A}(E, Y)$,
  $f$ induces $f_1$ and $f_0$ such that the above diagram commutes.
  Thus  $f_i = \Hom_{\A}(E, f_i^\prime)$ for some $f_i^\prime\in\Hom_\A(T_i, T_i^\prime)$,  $i = 0,
  1$.  So we get the following  diagram
  \[ \xymatrix{
  T_1 \ar[r]^u  \ar[d]^{f_1^\prime} &T_0 \ar[r]^\pi  \ar[d]^{f_0^\prime}
  &X
  \ar[r] \ar@{-->}[d]^{f^\prime} &0\\
T_1^\prime \ar[r]^{u^\prime}   &T_0^\prime \ar[r]^{\pi^\prime} &Y
  \ar[r]  &0\\
  }\]
with commutative left square.    So there exists
$f^\prime\in\Hom_\A(X, Y)$ such that the above diagram commutes.
Thus $f\Hom_{\A}(E,
  \pi) = \Hom_{\A}(E, f^\prime)\Hom_{\A}(E, \pi)$ and hence $f = \Hom_{\A}(E,
f^\prime)$. This proves that $\Hom_{\A}(E, -)$ is full. $\s$

\vskip5pt

\subsection{} The following result shows that, after taking the opposite algebra of
endomorphism algebra of a Gorenstein-projective generator, the
category of the Gorenstein-projective modules can  not be
``enlarged".

\vskip5pt

\begin{thm} \label{exp} \ Let $\A$ be an Artin algebra and $E$ a generator,
and $\Gamma = (\End_\A E)^{op}$. Then $\Hom_{\A}(E, -)$
induces an equivalence between $\mathcal{G}(E)$ and
$\mathcal{GP}(\Gamma)$.

\vskip5pt

In particular, if $E$ is a Gorenstein-projective generator, then
$\mathcal{GP}(\Gamma)$ is equivalent to a full subcategory of
$\mathcal{GP}(\A)$.
\end{thm}
\noindent {\bf Proof.} \ Let $L \in \mathcal{G}(E)$. By definition
there exists an exact complex $X^\bullet = (X^i, d^i)$ with all $X^i
\in \add E$ such that $X^\bullet$  is $E$-biexact and $L \cong \Ima
d^0$. Since $\Hom_{\A}(X^\bullet, E) \cong
\Hom_{\Gamma}(\Hom_{\A}(E, X^\bullet), \Gamma)$, it follows that
$\Hom_{\Gamma}(\Hom_{\A}(E, X^\bullet), \Gamma)$ is exact. This
means that  $\Hom_{\A}(E, X^\bullet)$ is a complete
$\Gamma$-projective resolution of $\Ima(\Hom_{\A}(E, d^0))$.
However, since $\Hom_{\A}(E, X^\bullet)$ is exact and E is a
generator, it follows that $\Ima(\Hom_{\A}(E, d^0))$ $ =
\Hom_{\A}(E, \Ima d^0)$. That is, $\Hom_{\A}(E, L) \in
\mathcal{GP}(\Gamma)$.

We claim that $\Hom_{\A}(E, -): \mathcal{G}(E) \longrightarrow
\mathcal{GP}(\Gamma)$ is dense. Suppose $Y \in
\mathcal{GP}(\Gamma)$. Then there is a complete projective
$\Gamma$-resolution $P^\bullet: \cdots \longrightarrow P^{-1}
\xlongrightarrow{d^{-1}} P^0 \xlongrightarrow {d^0}
P^{1}\xlongrightarrow{d^{1}} \cdots$ such that $Y \cong \Ima d^0.$
Then there exists a complex $E^\bullet : \cdots \longrightarrow
E^{-1} \xlongrightarrow{e^{-1}} E^0 \xlongrightarrow{e^0} E^{1}
\xlongrightarrow{e^1} \cdots$ with all $E^i \in \add E$ such that
$\Hom_{\A}(E, E^\bullet) = P^\bullet$. Since $E$ is a generator, for
each $i$ there is a surjective $\A$-map $\pi: E'\twoheadrightarrow
\Ker e^{i+1}$ with $E'\in\add E$; while $\Hom_{\A}(E', E^\bullet)$
is exact, so there is $g\in\Hom_{\A}(E', E^i)$ such that $\pi =
e^ig$. Decomposing $e^i$ as the composition $E^i \twoheadrightarrow
\Ima e^i \hookrightarrow \Ker e^{i+1}$, we see that the embedding
$\Ima e^i \hookrightarrow \Ker e^{i+1}$ is surjective. This proves
that $E^\bullet$ is also exact. Since $\Hom_{\Gamma}(P^\bullet,
\Gamma) \cong \Hom_{\A}(E^\bullet, E)$ and $\Hom_{\Gamma}(\Gamma,
P^\bullet) \cong \Hom_{\A}(E, E^\bullet)$, it follows that
$E^\bullet$ is $E$-biexact, and hence $\Ima e^i \in \mathcal{G}(E)$
for all $i$. Since $\Hom_{\A}(E, E^\bullet)$ is exact, $\Ima d^i
\cong \Ima (\Hom_{\A}(E, e^i)) = \Hom_{\A}(E, \Ima e^i)$. Thus $Y
\cong \Ima d^0 \in \Hom_{\A}(E, \mathcal{G}(E))$. This proves the
claim.

By Lemma \ref{yoneda} $\Hom_{\A}(E, -): \A\mbox{-}{\rm mod}
\longrightarrow  \Gamma\mbox{-}{\rm mod}$ is fully faithful. Thus
$\Hom_{\A}(E, -): \mathcal{G}(E) \longrightarrow
\mathcal{GP}(\Gamma)$ is an equivalence of categories.

\vskip5pt

If $E$ is a Gorenstein-projective generator, then
$\mathcal{G}(E)\subseteq \mathcal{GP}(\A)$ by Corollary \ref{GPE},
and hence $\mathcal{GP}(\Gamma)$ is equivalent to a full subcategory
of $\mathcal{GP}(\A)$. $\s$

\vskip5pt
\begin{cor} \ Let $\A$ be a {\rm CM}-finite Artin algebra, $E$ a Gorenstein-projective generator.
Then  $\Gamma = (\End E)^{op}$ is also {\rm CM}-finite.
\end{cor}
\subsection {} Recall that an Artin algebra is {\it basic}, if
$_\A\A$ is a direct sum of pairwise non-isomorphic indecomposable
modules. Any Artin algebra is Morita equivalent to a basic Artin
algebra.

\vskip5pt

Now we can prove the main result of this section, where the
assertion $(ii)$ is well-known.

\begin{thm} \label{free} \ Let $\A$ be a {\rm CM}-finite algebra, $E_1, \dots, E_n$ all the pairwise
non-isomorphic indecomposable Gorenstein-projective $\A$-modules.
Put $E = \bigoplus\limits_{1\le i\le n} E_i$ and ${\rm \bf Aus}
(\A)$ $ = (\End E)^{op}$. Then ${\rm \bf Aus} (\A)$ is {\rm
CM}-free.

\vskip5pt

More precisely, put $\Omega$ to be the class of pairwise
non-isomorphic basic {\rm CM}-finite algebras, and $\Theta$ the
class of pairwise non-isomorphic basic {\rm CM}-free algebras. Then

 \vskip5pt

$(i)$ \ The map ${\rm \bf Aus}: \ \Omega \longrightarrow \Theta$ is
surjective.

\vskip5pt

$(ii)$ \ The map ${\rm \bf Aus}$ sends  {\rm CM}-finite Gorenstein
algebras to algebras of finite global dimension.

\vskip5pt

$(iii)$ \ The map ${\rm \bf Aus}$ sends {\rm CM}-finite
non-Gorenstein algebras to  {\rm CM}-free  non-Gorenstein algebras.

\vskip5pt

$(iv)$ \ A {\rm CM}-free Gorenstein algebra $\A$ is of finite global
dimension.
\end{thm}
\noindent {\bf Proof.} Write $\Gamma = {\rm \bf Aus}(\A)$. By
definition we have $\mathcal{GP}(\A) = \add E \subseteq
\mathcal{G}(E)$. By Corollary \ref{GPE} $\mathcal{G}(E) \subseteq
\mathcal{GP}(\A)$. Thus $\mathcal{G}(E) = \add E$. By Theorem
\ref{exp} we have
$$\mathcal{GP}(\Gamma) =
\Hom_\A(E, \ \mathcal{G}(E)) = \Hom_\A(E, \ \add E) =
\mathcal{P}(\Gamma).$$ This means that $\Gamma$ is CM-free.

\vskip5pt

$(i)$ \ The map ${\rm \bf Aus}$ sends a basic CM-free algebra to
itself, and from this $(i)$ follows.

\vskip5pt

$(ii)$ \ Recall a well-known fact: for a CM-finite algebra $\A$,
$\A$ is Gorenstein if and only if $\Gamma$ is of finite global
dimension (see [LZ]). From this the assertion $(ii)$ follows.

\vskip5pt

$(iii)$ \ By the fact stated above, ${\rm \bf Aus}$ sends CM-finite
non-Gorenstein algebra $\A$ to a CM-free algebra $\Gamma$ of
infinite global dimension. Note that $\Gamma$ can not be Gorenstein:
otherwise, by $(ii)$ ${\rm \bf Aus}(\Gamma) = \Gamma$ is of finite
global dimension. This proves $(iii)$.

\vskip5pt

$(iv)$ \ By $(ii)$  ${\rm \bf Aus}(\A)$ is of finite global
dimension. Since $\A$ is Morita equivalent to ${\rm \bf Aus}(\A)$,
$\A$  is of finite global dimension. $\s$

\vskip5pt

\section{\bf Description of bounded derived categories}

Throughout this section $\mathcal A$ is an abelian category with
enough projective objects unless stated othwise. The bounded derived
category $D^b(\mathcal A)$ can be interpreted as $K^{-, b}(\mathcal
P(\mathcal A))$. The aim of this section is to describe
$D^b(\mathcal A)$ via ${\mathcal{GP}(A)}$, under the condition that
$\mathcal{A}$ is {\rm CM}-contravariantly finite.

\vskip5pt

\subsection{} To
describe $D^b(\mathcal A)$ via ${\mathcal{GP}(A)}$, we need the full
additive subcategory $K^{-, gpb}(\mathcal{GP}(\mathcal A))$ of
$K^{-}(\mathcal{GP}(\mathcal A))$ ([GZ, 3.3]), where
\begin{align*}  K^{-, gpb}(\mathcal{GP}(\mathcal A)) =   \{ \ & X^\bullet \in K^{-}(\mathcal{GP}(\mathcal A))
\mid \mbox{there exists an integer} \ N  \ \mbox{such that}  \\
& {\rm H}^{n}\Hom_\mathcal A(G, X^\bullet) = 0, \ \forall \ n\le N,
\ \forall \ G \in \mathcal{GP}(\mathcal A)\}. \end{align*}

\noindent As pointed in \cite{GZ}, $K^{-, gpb}(\mathcal{GP}(\mathcal
A))$ is a triangulated subcategory of $K^{-}(\mathcal{GP}(\mathcal
A))$. Clearly, it is the Gorenstein version of $K^{-,
b}(\mathcal{P}(\mathcal A))$.

\vskip5pt

Let $X^\bullet\in K^{-, gpb}(\mathcal{GP}(\mathcal A))$. Since
$\mathcal A$ has  enough projective objects and
$\mathcal{P}(\mathcal A) \subseteq \mathcal{GP}(\mathcal A)$, we see
that ${\rm H}^n(X^\bullet) = 0$ for $n< < 0$.

\vskip5pt

The main result of this section is as follows.

\vskip5pt

\begin{thm} \label{DB} \ Let $\mathcal{A}$ be a
{\rm CM}-contravariantly finite abelian category. Then there is a
triangle-equivalence $$D^b(\mathcal A)\cong {K^{-,
gpb}(\mathcal{GP}(\mathcal A))}/{K^{b, ac}(\mathcal{GP}(\mathcal
A))}$$ which induces a triangle-equivalence $\langle
\mathcal{GP}(\mathcal A)\rangle \cong K^{b} (\mathcal{GP}(\mathcal
A))/ K^{b, ac}(\mathcal{GP}(\mathcal A)).$
\end{thm}

\vskip5pt

\subsection{} To prove Theorem \ref{DB}, we need some preparations.
First, we have

\vskip5pt

\begin{lem} \label{b} $(i)$ \ ${K^{b, ac}(\mathcal{GP}(\mathcal A))}$ is a triangulated subcategory of $K^{-, gpb}(\mathcal{GP}(\mathcal A))$.

\vskip5pt

$(ii)$ \  Let $G^\bullet= (G^i, d^i) \in K^{-,
gpb}(\mathcal{GP}(\mathcal A))$. If $G^\bullet$ is exact, then
$G^\bullet \in {K^{b, ac}(\mathcal{GP}(\mathcal A))}$.
\end{lem}
\noindent {Proof.} \ $(i)$ is clear. We prove $(ii)$. Since
$\mathcal{GP}(\mathcal A)$ is closed under kernels of epimorphisms,
$\Ima{d^i} \in \mathcal{GP}(\mathcal A), \ \forall \ i\in \Bbb Z$.
By definition there exists an integer $N$ such that ${\rm
H}^{n}\Hom_{\mathcal A}(G, G^\bullet) = 0, \ \forall \ n \le N, \
\forall \ G \in \mathcal{GP}(\mathcal A)$. In particular ${\rm
H}^{n}\Hom_{\mathcal A}(\Ima{d^{n-1}}, G^\bullet) = 0.$ This implies
that the induced epimorphism $\widetilde{d^{n-1}}: G^{n-1}
\longrightarrow \Ima d^{n-1}$ splits for $n \le N$, and hence
$G^\bullet\cong G'^\bullet \in {K^{b, ac}(\mathcal{GP}(\mathcal
A))}$, where $G'^\bullet$ is the complex $\cdots \longrightarrow 0
\longrightarrow \Ima d^{N-1} \hookrightarrow G^{N} \longrightarrow
G^{N+1}\longrightarrow \cdots. $
 $\s$

\vskip5pt

The following observation plays an important role in our argument.

\vskip5pt

\begin{lem} \label{quasi} Assume
that  $\mathcal A$ is {\rm CM}-contravariantly finite  abelian
category, $P^\bullet \in K^{-, b}(\mathcal{P}(\mathcal A))$.  Then

\vskip5pt

$(i)$ \ There exists a quasi-isomorphism $P^\bullet \ra G^\bullet$
with $G^\bullet \in K^{-, gpb}(\mathcal{GP}(\mathcal A))$.

\vskip5pt

$(ii)$ \ For chain maps $f_i^\bullet: P^\bullet \ra G_i^\bullet$
with $G_i^\bullet \in K^{-, gpb}(\mathcal{GP}(\mathcal A)), \ 1\le
i\le m$, there exist a quasi-isomorphism $g^\bullet: P^\bullet \ra
G^\bullet$ and chain maps $h_i^\bullet: G^\bullet  \ra G^\bullet_i$
with $G^\bullet \in K^{-, gpb}(\mathcal{GP}(\mathcal A))$, such that
$f_i^\bullet  = h_i^\bullet g^\bullet, \ 1\le i\le m$.
\end{lem}

\noindent{\bf Proof.} \ $(i)$ \ Write $P^\bullet = (P^i, d^i)$. Let
$N$ be an integer such that ${\rm H}^{n} P^\bullet = 0$ for $n \le
N$. Since $\mathcal{GP}(\mathcal A)$ is  a contravariantly finite
subcategory in $\mathcal A$, we can  take a right
$\mathcal{GP}(\mathcal A)$-approximation $G^{N-1} \longrightarrow
{\rm Ker}d^{N}$ of  ${\rm Ker}d^{N}$, and then take right
$\mathcal{GP}(\mathcal A)$-approximations  step by step
\[
\xymatrix @C=0.5pc {\cdots & & \ar[rr] & & G
^{N-2}\ar[rd]\ar[rr]^{\partial^{N-2}} & & G^{N-1}\ar[rd] \ar[rr] & &
P^{N} \\ & && {\rm Ker}\partial^{N-2}\ar[ru] & & {\rm
Ker}\partial^{N-1} \ar[ru] & & {\rm Ker}d^{N}\ar[ru].}
\] In this way we get a complex $$G^\bullet:  \cdots \ra G^{N-2} \ra G^{N-1} \ra P^{N} \ra
P^{N+1} \ra \cdots$$ in $K^{-, gpb}(\mathcal{GP}(\mathcal A))$ with
${\rm H}^{n}\Hom_{\mathcal A}(G, G^\bullet) = 0, \forall \ n \le N,
\ \forall \ G \in \mathcal{GP}(\mathcal A)$. Since $G^{n}$ with
$n\le N$ is constructed via right $\mathcal{GP}(\mathcal
A)$-approximations, it is clear that there exists an induced a chain
map
$$\xymatrix{P^\bullet= \ar@<-1.2ex>@{-->}[d] &\cdots \ar[r] &P^{N-2} \ar@{-->}[d] \ar[r] &P^{N-1} \ar@{-->}[d] \ar[r] &P^{N} \ar@{=}[d]  \ar[r] &P^{N+1}
\ar@{=}[d]  \ar[r] &\cdots\\
G^\bullet= &\cdots \ar[r] &G^{N-2} \ar[r]&G^{N-1} \ar[r] &P^{N}
\ar[r]  &P^{N+1} \ar[r]  &\cdots,}$$ which is a quasi-isomorphism,
since ${\rm H}^{n} P^\bullet = 0 = {\rm H}^{n}G^\bullet$ for $n \le
N$.

\vskip5pt

$(ii)$ \ Let $N$ be an integer such that ${\rm H}^{n}\Hom_{\mathcal
A}(G, G^\bullet_i) = 0 = {\rm H}^{n}P^\bullet, \ \forall \ n \le N,
\ \forall \ G \in \mathcal{GP}(\mathcal A), \ 1\le i\le m$. By the
proof of $(i)$ we have the following quasi-isomorphism
$$\xymatrix{P^\bullet = \ar@<-1.2ex>[d]^{g^\bullet} &\cdots \ar[r] &P^{N-2} \ar[d] \ar[r] &P^{N-1} \ar[d] \ar[r] &P^{N} \ar@{=}[d]  \ar[r] &P^{N+1} \ar@{=}[d]  \ar[r] &\cdots\\
G^\bullet =  &\cdots \ar[r] &G^{N-2} \ar[r]&G^{N-1} \ar[r] &P^{N}
\ar[r]  &P^{N+1} \ar[r]  &\cdots}$$ with $G^\bullet \in K^{-,
gpb}(\mathcal{GP}(\mathcal A))$. Since for each $i$ we have ${\rm
H}^{n}\Hom_{\mathcal A}(G, G^\bullet_i) = 0, \forall \ n \le N, \
\forall \ G \in \mathcal{GP}(\mathcal A)$, by applying
$\Hom_{\mathcal A}(G^{N-1}, -), \ \Hom_{\mathcal A}(G^{N-2}, -), \
\cdots,$ \ to $G_i^\bullet$ respectively, we obtain the following
chain map $h_i^\bullet$:
$$\xymatrix{G^\bullet = \ar@<-1.2ex>@{-->}[d]^{h_i^\bullet} &\cdots \ar[r] &G^{N-2} \ar@{-->}[d] ^{h_i^{N-2}}
\ar[r] &G^{N-1} \ar@{-->}[d] ^{h_i^{N-1}}  \ar[r]  &P^{-N}
\ar[d]^{f_i^{N}}
  \ar[r] &P^{N+1} \ar[d]^{f_i^{N+1}}  \ar[r] &\cdots\\
G^\bullet_i = &\cdots \ar[r] &G_i^{N-2} \ar[r]&G_i^{N-1} \ar[r]
&G_i^{N} \ar[r]  &G_i^{N+1} \ar[r]  &\cdots.}$$ For  $1\le i\le m$,
if $l \ge N$, then $f_i^{l}-h_i^{l} g^{l}=0$. If $l = N-1$, then
$f_i^{N-1}-h_i^{N-1} g^{N-1}$ can  factor through $G_i^{N-2}
\longrightarrow G_i^{N-1}$ since ${\rm H}^{N-1} G^\bullet_i = 0$. By
induction we get the following null homotopy
$f_i^\bullet-h_i^\bullet g^\bullet$ .
$$\xymatrix{P^\bullet = \ar@<-1.2ex>[d]^{f_i^\bullet-h_i^\bullet g^\bullet} &\cdots \ar[r] &P^{N-2} \ar@{-->}[dl]  \ar[d]
\ar[r] &P^{N-1} \ar@{-->}[dl]  \ar[d]   \ar[r]  &P^{N}
\ar@{-->}[dl]_{0} \ar[d]^{0}
  \ar[r] &P^{N+1} \ar@{-->}[dl]_{0} \ar[d]^{0}  \ar[r] &\cdots\\
G^\bullet_i = &\cdots \ar[r] &G_i^{N-2} \ar[r]&G_i^{N-1} \ar[r]
&G_i^{N} \ar[r]  &G_i^{N+1} \ar[r]  &\cdots.}$$ This completes the
proof. $\s$

\vskip5pt

\subsection{} {\bf Proof of Theorem \ref{DB}.} Let $\eta: K^{-, gpb}(\mathcal{GP}(\mathcal
A)) \longrightarrow D^-(\mathcal A)$ be the composition
$$K^{-, gpb}(\mathcal{GP}(\mathcal A)) \stackrel {\sigma} \hookrightarrow  K^{-}(\mathcal A) \xlongrightarrow {Q} D^-(\mathcal A) = K^{-}(\mathcal A)/K^{-, ac}(\mathcal A)$$
where $\sigma$ is the embedding and $Q$ is the localization functor.

Since $K^{b, ac}(\mathcal{GP}(\mathcal A))$ is a triangulated
subcategory of $K^{-, gpb}(\mathcal{GP}(\mathcal A))$, we have the
Verdier quotient ${ K^{-, gpb}(\mathcal{GP}(\mathcal A))}/{K^{b,
ac}(\mathcal{GP}(\mathcal A))}$. Since $\eta(K^{b,
ac}(\mathcal{GP}(\mathcal A))) = 0$, by the universal property
$\eta$ induces a unique triangle functor $\overline{\eta}: K^{-,
gpb}(\mathcal{GP}(\mathcal A))/K^{b, ac}(\mathcal{GP}(\mathcal
A))\longrightarrow D^-(\mathcal A)$.

\vskip5pt

If $\mathcal A$ is CM-contravariantly finite, then by Lemma
\ref{quasi}$(i)$ we have ${\rm Im}\overline{\eta} = K^{-,
b}(\mathcal{P}(\mathcal A))\cong D^b(\mathcal A)$. So we get a dense
triangle functor from $K^{-, gpb}(\mathcal{GP}(\mathcal A))/K^{b,
ac}(\mathcal{GP}(\mathcal A))$ to $D^b(\mathcal A)$, again denoted
by $\overline{\eta}$.

\vskip5pt

We prove that $\overline{\eta}$ is fully faithful. Suppose
$G^\bullet \in K^{-, gpb}(\mathcal{GP}(\mathcal A))$ and $\overline
\eta(G^\bullet) = 0$. Then $G^\bullet$ is exact.  By Lemma \ref{b}
we have $G^\bullet \in {K^{b, ac}(\mathcal{GP}(\mathcal A))}$. That
is $\overline{\eta}$ maps non-zero objects to non-zero objects. By
Lemma \ref{ff} it suffices to prove that $\overline{\eta}$ is full.

\vskip5pt

Let $G^\bullet_1, G^\bullet_2 \in K^{-, gpb}(\mathcal{GP}(\mathcal
A))$, and ${\alpha^\bullet}/{s^\bullet}$ be a morphism in
$\Hom_{D^-(\mathcal A)}(\overline{\eta}(G^\bullet_1),
\overline{\eta}(G^\bullet_2)) = \Hom_{D^-(\mathcal A)}(G^\bullet_1,
G^\bullet_2)$, where $s^\bullet: X^\bullet \Longrightarrow
G^\bullet_1$ is a quasi-isomorphism with $X^\bullet\in K^-(\mathcal
A)$, and $\alpha^\bullet: X^\bullet \longrightarrow  G^\bullet_2$ is
a morphism in $K^-(\mathcal A)$. Then there exists a
quasi-isomorphism $t^\bullet: P^\bullet \Longrightarrow X^\bullet$
with $P^\bullet \in K^{-}(\mathcal{P}(\mathcal A))$. Since
$s^\bullet$ and $t^\bullet$ are quasi-isomorphisms and $G^\bullet_1
\in K^{-, gpb}(\mathcal{GP}(\mathcal A))$, it follows that
$P^\bullet\in K^{-, b}(\mathcal{P}(\mathcal A))$. Thus we get the
commutative diagram
$$\xymatrix{&  X^\bullet \ar@{=>}[dl]_{s^\bullet} \ar[dr]^{\alpha^\bullet}\\ G^\bullet_1  &P^\bullet \ar@{=>}[l]_-{s^\bullet t^\bullet} \ar[r]^-{\alpha^\bullet t^\bullet} \ar@{=>}[u]^-{t^\bullet}
&G^\bullet_2}$$ where the double arrowed morphisms mean
quasi-isomorphisms. By Lemma \ref{quasi}$(ii)$ we have the following
commutative diagram
 $$\xymatrix{ G^\bullet_1  &P^\bullet \ar@{=>}[l]_- {s^\bullet t^\bullet} \ar[r]^-{\alpha^\bullet t^\bullet} \ar@{=>}[d]^{g^\bullet}
 &G^\bullet_2 \\ & G^\bullet \ar@{=>}[ul]^{l^\bullet} \ar[ur]_{\beta^\bullet}}$$
where $G^\bullet\in K^{-, gpb}(\mathcal{GP}(\mathcal A))$, and
$g^\bullet: P^\bullet\Longrightarrow G^\bullet$ is a
quasi-isomorphism. Note that $l^\bullet$ is also a
quasi-isomorphism, hence the mapping cone ${\rm Cone}(l^\bullet)$ is
exact. Since $K^{-, gpb}(\mathcal{GP}(\mathcal A))$ is a
triangulated subcategory of $K^{-}(\mathcal A)$, ${\rm
Cone}(l^\bullet) \in K^{-, gpb}(\mathcal{GP}(\mathcal A))$. By Lemma
\ref{b}$(ii)$ ${\rm Cone}(l^\bullet) \in K^{b,
ac}(\mathcal{GP}(\mathcal A))$. This proves
${\beta^\bullet}/{l^\bullet}\in \Hom_{K^{-,
gpb}(\mathcal{GP}(\mathcal A))/K^{b, ac}(\mathcal{GP}(\mathcal
A))}(G^\bullet_1, G^\bullet_2)$ and ${\alpha^\bullet}/{s^\bullet} =
{\beta^\bullet}/{l^\bullet} =
\overline{\eta}({\beta^\bullet}/{l^\bullet})$. This proves the first
triangle-equivalence in Theorem \ref{DB}.

\vskip5pt

By Lemma \ref{verdierquotient} we know that
$K^{b}(\mathcal{GP}(\mathcal A))/ K^{b, ac}(\mathcal{GP}(\mathcal
A))$ is a triangulated subcategory of ${K^{-,
gpb}(\mathcal{GP}(\mathcal A))}/{K^{b, ac}(\mathcal{GP}(\mathcal
A))}$. Thus $\overline{\eta}(K^{b}(\mathcal{GP}(\mathcal A))/ K^{b,
ac}(\mathcal{GP}(\mathcal A)))$ is a triangulated subcategory of
$D^b(\mathcal A)$. Obviously this image is generated by
$\mathcal{GP}(\mathcal A)$, i.e.,
$$\overline{\eta}(K^{b}(\mathcal{GP}(\mathcal A))/ K^{b, ac}(\mathcal{GP}(\mathcal A)))
= \langle \mathcal{GP}(\mathcal A)\rangle.$$ So we get the second
triangle-equivalence in Theorem \ref{DB}. $\s$

\vskip5pt

\section{\bf Descriptions of Gorenstein defect categories}

Throughout this section $\mathcal A$ is an abelian category with
enough projective objects. The aim of this section is to describe
$D^b_{{\rm defect} }(\mathcal A)$.

\vskip5pt

\subsection{}  A quick description is

\vskip5pt

\begin{lem} \label{defect} \ Let $F:  \underline
{\mathcal{GP}(\mathcal A)} \longrightarrow D^b_{sg}(\mathcal A) : =
D^b(\mathcal A)/K^b(\mathcal P(\mathcal A))$ be the fully-faithful
triangle functor defined in $(2.1)$. Then $${\rm Im}F = \langle
\mathcal{GP}(\mathcal A)\rangle /K^b(\mathcal P(\mathcal A)),$$ and
hence we have a triangle-equivalence
$$D^b_{{\rm defect}}(\mathcal A) \cong D^b(\mathcal A)/\langle \mathcal{GP}(\mathcal A)\rangle.$$
\end{lem} \noindent{\bf Proof.} \ Since $\langle
\mathcal{GP}(\mathcal A)\rangle $ is the triangulated subcategory of
$D^b(\mathcal A)$ generated by $\mathcal{GP}(\mathcal A)$,  it
follows that $\langle \mathcal{GP}(\mathcal A)\rangle/K^b(\mathcal
P(\mathcal A))$ is the triangulated subcategory of $D^b(\mathcal
A)/K^b(\mathcal P(\mathcal A))$ generated by $\mathcal{GP}(\mathcal
A)$, here view objects in $\mathcal{GP}(\mathcal A)$ as  stalk
complexes at degree $0$. Since $F$ is full, it follows that ${\rm
Im}F$ is a triangulated subcategory of $D^b(\mathcal A)/K^b(\mathcal
P(\mathcal A))$ containing $\mathcal{GP}(\mathcal A)$. It follows
that ${\rm Im}F \supseteq \langle \mathcal{GP}(\mathcal
A)\rangle/K^b(\mathcal P(\mathcal A)) .$

\vskip5pt

On the other hand, by definition ${\rm Im}F \subseteq\langle
\mathcal{GP}(\mathcal A)\rangle/K^b(\mathcal P(\mathcal A)).$ It
follows that ${\rm Im}F = \langle \mathcal{GP}(\mathcal A)\rangle
/K^b(\mathcal P(\mathcal A)),$ and hence by Lemma
\ref{verdierquotient} we have a triangle-equivalence $$ \ \
D^b_{{\rm defect}}(\mathcal A) = (D^b(\mathcal A)/K^b(\mathcal
P(\mathcal A)))/(\langle \mathcal{GP}(\mathcal A)\rangle
/K^b(\mathcal P(\mathcal A))) \cong D^b(\mathcal A)/\langle
\mathcal{GP}(\mathcal A)\rangle. \ \ \ \ \s $$

\vskip5pt

By definition we have $K^{b} (\mathcal{P}(\mathcal A))\subseteq
\langle \mathcal{GP}(\mathcal A)\rangle\subseteq D^b(\mathcal A).$

\vskip5pt

\begin{cor} \label{cordefect} \ We have

\vskip5pt

$(i)$ \ $K^{b} (\mathcal{P}(\mathcal A)) = \langle
\mathcal{GP}(\mathcal A)\rangle$ if and only if $\mathcal A$ is {\rm
CM}-free, i.e., $\mathcal{P}(\mathcal A) = \mathcal{GP}(\mathcal
A)$.

\vskip5pt

$(ii)$ \ $\langle \mathcal{GP}(\mathcal A)\rangle = D^b(\mathcal A)$
if and only if each object of $\mathcal A$ has finite
Gorenstein-projective dimension.
\end{cor}
\noindent{\bf Proof.} $(i)$ If $\mathcal{P}(\mathcal A) =
\mathcal{GP}(\mathcal A)$, then it is clear that $K^{b}
(\mathcal{P}(\mathcal A)) = \langle \mathcal{GP}(\mathcal
A)\rangle$. Conversely, assume that $K^{b} (\mathcal{P}(\mathcal A))
= \langle \mathcal{GP}(\mathcal A)\rangle$. Then one easily see that
any Gorenstein-projective object $G$ is of finite projective
dimension, and hence $G$ is projective ([EJ2, 10.2.3]).

\vskip5pt

$(ii)$ By Lemma \ref{defect} $\langle \mathcal{GP}(\mathcal
A)\rangle = D^b(\mathcal A)$ if and only if $D^b_{{\rm
defect}}(\mathcal A) = 0$, and if and only if each object of
$\mathcal A$ has finite Gorenstein-projective dimension. $\s$

\vskip5pt

\begin{rem} Although the generating process inside a triangulated category
is clear {\rm (see e.g. R. Rouquier [Rou])}, but the output of this
generating process is not so clear. In other words, the description
in Lemma \ref{defect} is rough in the sense that the shape of
$\langle \mathcal{GP}(\mathcal A)\rangle$ is vague. We will give an
explicit description of $D^b_{{\rm defect} }(\mathcal A)$.
\end{rem}

\subsection{} Define $K^{-, \ b}_{\mathcal G}(\mathcal
P(\mathcal A))$ to be the full subcategory of $K^{-, \ b}(\mathcal
P(\mathcal A))$ consisting of all the complexes $(P^\bullet, d)\in
K^{-, \ b}(\mathcal P(\mathcal A))$ such that there exists
$n_0\in\Bbb Z$ with ${\rm H}^m(P^\bullet ) = 0, \ \forall \ m\le
n_0$, and ${\rm Ker} d^{n_0}\in \mathcal{GP}(\mathcal A).$

\vskip5pt

In order to say that $K^{-, \ b}_{\mathcal G}(\mathcal P(\mathcal
A))$ is a thick triangulated subcategory of $K^{-, \ b} (\mathcal
P(\mathcal A))$, first, we need to say that $K^{-, \ b}_{\mathcal
G}(\mathcal P(\mathcal A))$ is closed under isomorphisms of $K^{-, \
b}(\mathcal P(\mathcal A))$.

\vskip5pt

The following fact is well-known.

\vskip5pt

\begin{lem} \label {isoinstable}  Let  $\underline {\mathcal A}$ be the stable
category of $\mathcal A$ modulo $\mathcal P(\mathcal A)$. Let $X,
Y\in \mathcal A$. Then $X\cong Y$ in $\underline {\mathcal A}$ if
and only if there exist projective objects $P$ and $Q$ such that
$X\oplus P \cong Y\oplus Q$ in $\mathcal A$.
\end{lem}

\vskip5pt

The following fact is also well-known.

\vskip5pt

\begin{lem}
\label{factor} Let $f^\bullet: P^\bullet \longrightarrow Q^\bullet$
be a null-homotopy with  $P^\bullet, \ Q^\bullet \in K^{-, \
b}(\mathcal P(\mathcal A))$.  Assume that ${\rm H}^m(P^\bullet ) = 0
= {\rm H}^m(Q^\bullet ), \ \forall \ m\le n_0$. Then the restriction
$\widetilde{f^{n_0}}: {\rm Ker}d^{n_0}_{P^\bullet}\longrightarrow
{\rm Ker}d^{n_0}_{Q^\bullet}$ of $f^{n_0}$ factors through a
projective object.
\end{lem}
\noindent {\bf Proof.} For convenience we include a proof.  Let $s:
f^\bullet \sim 0$ be a homotopy. Put $K ={\rm
Ker}d^{n_0}_{P^\bullet}$ and ~$L ={\rm Ker}d^{n_0}_{Q^\bullet}$.

$${\xymatrix{\cdots\ar[rr] && P^{{n_0}-1}\ar[rr]\ar[dd]^{f^{{n_0}-1}}
\ar[dr]^{\pi}&& P^{n_0}\ar[rr]^{d^{n_0}}\ar[dd]^{f^{{{n_0}}}}\ar@/^/@{.>}[ddll]^{\quad s^{n_0}} && P^{{n_0}+1}\ar[rr]\ar@{.>}[ddll]^{s^{{n_0}+1}} && \cdots\\
&&&K\ar[ur]^i\ar[dd]\\
\cdots\ar[rr] && Q^{{n_0}-1}\ar[rr]\ar[dr]^{\pi^\prime} && Q^{n_0}\ar[rr]^{d^{n_0}} && Q^{{n_0}+1}\ar[rr] && \cdots\\
&&&L\ar[ur]^{i^\prime}}}$$

\noindent Then we have ~$f^{n_0} i=i^\prime \widetilde{f^{n_0}}$ and
~$\widetilde{f^{n_0}}\pi=\pi^\prime f^{n_0-1}$. Since
$$
i^\prime\widetilde{f^{n_0}} \pi=d^{n_0-1}f^{n_0-1}=d^{n_0-1}(s^{n_0}
d^{n_0-1} + d^{n_0-2} s^{n_0-1})=d^{n_0-1} s^{n_0} d^{n_0-1}=
i^\prime \pi^\prime s^{n_0} i \pi,$$ we have
~$\widetilde{f^{n_0}}=\pi^\prime s^{n_0} i$. Thus
$\widetilde{f^{n_0}}$ factors through a projective object. $\s$

\vskip5pt

The following fact shows that $K^{-, \ b}_{\mathcal G}(\mathcal
P(\mathcal A))$ is closed under isomorphisms of $K^{-, \ b}(\mathcal
P(\mathcal A))$.

\vskip5pt

\begin{lem} \label{basic} Let $f^\bullet: P^\bullet \longrightarrow Q^\bullet$
be a homotopy equivalence in $K^{-, \ b}(\mathcal P(\mathcal A))$
with $P^\bullet \in K^{-, \ b}_{\mathcal G}(\mathcal P(\mathcal
A))$. Then $Q^\bullet \in K^{-, \ b}_{\mathcal G}(\mathcal
P(\mathcal A))$.
\end{lem}
\noindent {\bf Proof.} Let $n_0\in\Bbb Z$ such that ${\rm
H}^m(P^\bullet ) = 0, \ \forall \ m\le n_0$, and ${\rm Ker}
d^{n_0}\in \mathcal{GP}(\mathcal A).$ Then ${\rm H}^m(Q^\bullet ) =
0, \ \forall \ m\le n_0$. Let  $g^\bullet: Q^\bullet\longrightarrow
P^\bullet$ be the inverse of $f^\bullet$. Then we have the
restriction $\widetilde{f^{n_0}}: {\rm
Ker}d^{n_0}_{P^\bullet}\longrightarrow {\rm Ker}d^{n_0}_{Q^\bullet}$
of $f^{n_0}$, and the restriction $\widetilde{g^{n_0}}: {\rm
Ker}d^{n_0}_{Q^\bullet}\longrightarrow {\rm Ker}d^{n_0}_{P^\bullet}$
of $g^{n_0}$. Since $g^\bullet f^\bullet-{\rm Id}_{P^\bullet}$ and
~$f^\bullet g^\bullet -{\rm Id}_{Q^\bullet}$ are null-homotopy, by
Lemma \ref{factor} the induced morphism
$\widetilde{g^{n_0}}\widetilde{f^{n_0}}-{\rm Id}_{{\rm
Ker}d^{n_0}_{P^\bullet}}$ factors through a projective object,
and~$\widetilde{f^{n_0}}\widetilde{g^{n_0}}-{\rm Id}_{{\rm
Ker}d^{n_0}_{Q^\bullet}}$ factors through a projective object. This
means that ${\rm Ker}d^{n_0}_{P^\bullet}$ and ${\rm
Ker}d^{n_0}_{Q^\bullet}$ are isomorphic in $\underline {\mathcal
A}$. By Lemma \ref{isoinstable} there exist projective objects $P$
and $Q$ such that there is an isomorphism ${\rm
Ker}d^{n_0}_{P^\bullet}\oplus P \cong {\rm
Ker}d^{n_0}_{Q^\bullet}\oplus Q$ in $\mathcal A$. Since ${\rm Ker}
d^{n_0}_{P^\bullet}\in \mathcal{GP}(\mathcal A)$ and
$\mathcal{GP}(\mathcal A)$ is closed under direct summands, it
follows that ${\rm Ker} d^{n_0}_{Q^\bullet}\in \mathcal{GP}(\mathcal
A)$. This proves $Q^\bullet \in K^{-, \ b}_{\mathcal G}(\mathcal
P(\mathcal A))$. $\s$

\vskip5pt

\begin{prop} \label{Gversion} \ $K^{-, \ b}_{\mathcal G}(\mathcal
P(\mathcal A))$ is a thick triangulated subcategory of $K^{-, \ b}
(\mathcal P(\mathcal A))$.
\end{prop}
\noindent {\bf Proof.} \ It is clear that $K^{-, \ b}_{\mathcal
G}(\mathcal P(\mathcal A))$ is an additive category.  By Lemma
\ref{basic} $K^{-, \ b}_{\mathcal G}(\mathcal P(\mathcal A))$ is a
full subcategory of $K^{-, \ b}(\mathcal P(\mathcal A))$ closed
under isomorphisms. Since $\mathcal{GP}(\mathcal A)$ is closed under
direct summands, it follows that $K^{-, \ b}_\mathcal G(\mathcal
P(\mathcal A))$ is closed under direct summands. It is also clear
that $K^{-, \ b}_{\mathcal G}(\mathcal P(\mathcal A))$ is closed
under the shift functor $[1]$ and $[-1]$. Let $f^\bullet: P^\bullet
\longrightarrow Q^\bullet$ be a chain map with $P^\bullet,
Q^\bullet\in K^{-, \ b}_{\mathcal G}(\mathcal P(\mathcal A)).$ It
remains to prove that the mapping ${\rm Cone}(f^\bullet)\in K^{-, \
b}_{\mathcal G}(\mathcal P(\mathcal A)).$

Since $P^\bullet, Q^\bullet\in K^{-, \ b}_{\mathcal G}(\mathcal
P(\mathcal A)),$ there exists an integer $n_0$ such that ${\rm
H}^m(P^\bullet ) = 0 = {\rm H}^m(Q^\bullet)$ for $m\le n_0$, and
${\rm Ker} d_{P^\bullet}^{n_0}\in \mathcal{GP}(A), \ {\rm Ker}
d_{Q^\bullet}^{n_0}\in \mathcal{GP}(\mathcal A)$. It follows that
there are complete projective resolutions $X^\bullet$ and
$Y^\bullet$ such that $${\rm Ker} d_{X^\bullet}^{n_0} = {\rm Ker}
d_{P^\bullet}^{n_0}, \ \ {\rm Ker} d_{Y^\bullet}^{n_0} = {\rm Ker}
d_{Q^\bullet}^{n_0}.$$ Construct the following two complexes
$$P'^\bullet = \cdots \longrightarrow P^{n_0-2}\stackrel{d^{n_0-2}_{P^\bullet}}
\longrightarrow P^{n_0-1}
\stackrel{d^{n_0-1}_{P'^\bullet}}\longrightarrow
X^{n_0}\stackrel{d^{n_0}_{P^\bullet}}\longrightarrow X^{n_0+1}
\stackrel{d^{n_0+1}_{X^\bullet}}\longrightarrow
X^{n_0+2}\longrightarrow\cdots,$$ and
$$Q'^\bullet = \cdots \longrightarrow Q^{n_0-2}\stackrel{d^{n_0-2}_{Q^\bullet}}
\longrightarrow Q^{n_0-1}
\stackrel{d^{n_0-1}_{Q'^\bullet}}\longrightarrow
Y^{n_0}\stackrel{d^{n_0}_{Q^\bullet}}\longrightarrow Y^{n_0+1}
\stackrel{d^{n_0+1}_{Y^\bullet}}\longrightarrow
Y^{n_0+2}\longrightarrow\cdots,$$ where $d^{n_0}_{P'^\bullet}:
P^{n_0-1}\longrightarrow X^{n_0}$ is the composition of the
canonical morphisms
$$P^{n_0-1}\twoheadrightarrow {\rm Ker}
d_{P^\bullet}^{n_0}, \ \ {\rm Ker}
d_{X^\bullet}^{n_0}\hookrightarrow X^{n_0},$$ and
$d^{n_0}_{Q'^\bullet}: Q^{n_0}\longrightarrow Y^{n_0+1}$ is the
composition of the canonical morphisms
$$Q^{n_0-1}\twoheadrightarrow {\rm Ker}
d_{Q^\bullet}^{n_0}, \ \ {\rm Ker}
d_{Y^\bullet}^{n_0}\hookrightarrow Y^{n_0}.$$ Now $f^\bullet:
P^\bullet \longrightarrow Q^\bullet$ induces a morphism ${\rm Ker}
d_{P^\bullet}^{n_0}\longrightarrow {\rm Ker} d_{Q^\bullet}^{n_0}$,
namely a morphism ${\rm Ker} d_{X^\bullet}^{n_0}\longrightarrow {\rm
Ker} d_{Y^\bullet}^{n_0}$. Since projective objects are injective
objects in the category $\mathcal{GP}(\mathcal A)$, it follows that
the ${\rm Ker} d_{X^\bullet}^{n_0}\longrightarrow {\rm Ker}
d_{Y^\bullet}^{n_0}$ induces morphisms $f'^i: X^i \longrightarrow
Y^i$ for $i\ge n_0$, such that the diagram
\[\xymatrix {P'^\bullet\ar[d]_{f'^\bullet}: &
\cdots\ar[r] &
P^{n_0-2}\ar[d]_{f^{n_0-2}}\ar[r]^{d_{P^\bullet}^{n_0-2}} &
P^{n_0-1}\ar[d]_{f^{n_0-1}}\ar[r]^{d_{P'^\bullet}^{n_0-1}} &
X^{n_0}\ar[d]_{f'^{n_0}}\ar[r]^{d_{P^\bullet}^{n_0}} & X^{n_0+1}
\ar[d]_{f'^{n_0+1}} \ar[r]^{d_{X^\bullet}^{n_0+1}} & \cdots
\\
Q'^\bullet: & \cdots\ar[r] & Q^{n_0-2}\ar[r]^{d_{Q^\bullet}^{n_0-2}}
& Q^{n_0-1}\ar[r]^{d_{Q'^\bullet}^{n_0-1}} &
Y^{n_0}\ar[r]^{d_{Q^\bullet}^{n_0}} &
Y^{n_0+1}\ar[r]^{d_{Y^\bullet}^{n_0+1}} & \cdots.}\] commutes. By
construction $P'^\bullet$ and $Q'^\bullet$ are complete projective
resolutions. It is clear that ${\rm Cone}(f'^\bullet)$ is again a
complete projective resolution, and in particular ${\rm Ker} \
d^{n_0-3}_{{\rm Cone}(f'^\bullet)}$ is a Gorenstein-projective
object. However by construction we have
$$d^{m}_{{\rm Cone}(f^\bullet)} =\left (
\begin{smallmatrix}-d_{P^\bullet}^{m+1} & 0 \\
f^{m+1} & d^{m}_{Q^\bullet}  \end{smallmatrix} \right )= d^{m}_{{\rm
Cone}(f'^\bullet)}, \ \forall \ m\le n_0-3.$$ This proves that ${\rm
Cone}(f^\bullet)\in K^{-, \ b}_\mathcal G(\mathcal{P}(\mathcal A)).$
$\s$

\vskip5pt

\subsection{} We have the following description of $D^b_{{\rm defect} }(\mathcal A)$.

\vskip5pt

\begin{thm} \label{descriptionofdeft} \ Let $\mathcal A$ be an
abelian category with enough projective objects. Then

\vskip5pt

$(i)$ \ There is a triangle-equivalence
$$D^b_{{\rm defect}}(\mathcal A) \cong K^{-, \ b}(\mathcal P(\mathcal A))/K^{-, \ b}_\mathcal G(\mathcal P(\mathcal A)).$$

\vskip5pt

$(ii)$ If in addition $\mathcal A$ is {\rm CM}-contravariantly
finite, then there is a triangle-equivalence$$D^b_{{\rm
defect}}(\mathcal A) \cong K^{-, gpb}(\mathcal {GP} (\mathcal
A))/K^{b}(\mathcal {GP} (\mathcal A)).$$
\end{thm}

\noindent {\bf Proof.} \ $(i)$ \ We claim: the restriction of the
canonical triangle-equivalence $\rho: K^{-, b}(\mathcal P(\mathcal
A)) \longrightarrow D^b(\mathcal A)$ to $K^{-, b}_{\mathcal
G}(\mathcal P(\mathcal A))$ gives rise to a triangle-equivalence
$K^{-, b}_{\mathcal G}(\mathcal P(\mathcal A)) \longrightarrow
\langle \mathcal{GP}(\mathcal A)\rangle.$

\vskip5pt

In fact, let $P^\bullet \in K^{-, \ b}_{\mathcal G}(\mathcal
P(\mathcal A))$. By definition there is an integer $n_0\in\Bbb Z$
such that ${\rm H}^m(P^\bullet ) = 0, \ \forall \ m\le n_0$, and
${\rm Ker} d^{n_0}\in \mathcal{GP}(\mathcal A).$ Then there is a
quasi-isomorphism $P^\bullet \longrightarrow G^\bullet$, where
$G^\bullet $ is the bounded complex
$$0 \longrightarrow {\rm Ker} d^{n_0}\hookrightarrow P^{n_0} \stackrel
{d^{n_0}_{P^\bullet}}\longrightarrow P^{n_0+1} \longrightarrow
\cdots.
$$  It follows that $\rho(P^\bullet) =
P^\bullet\cong G^\bullet$ in $D^b(\mathcal A)$. Since $P^\bullet$ is
a upper-bounded complex of projective objects, it follows that $G$
is a bounded complex of Gorenstein-projective objects. Thus
$G^\bullet\in \langle \mathcal{GP}(\mathcal A)\rangle$. So $\rho
(K^{-, b}_{\mathcal G}(\mathcal P(\mathcal A))) \subseteq \langle
\mathcal{GP}(\mathcal A)\rangle.$

\vskip5pt

On the other hand, by Proposition \ref{Gversion}  $K^{-, \
b}_{\mathcal G}(\mathcal P(\mathcal A))$ is a triangulated
subcategory of $K^{-, \ b} (\mathcal P(\mathcal A))$. Since $\rho:
K^{-, b}(\mathcal P(\mathcal A)) \longrightarrow D^b(\mathcal A)$ is
full, it follows that $\rho (K^{-, b}_{\mathcal G}(\mathcal
P(\mathcal A)))$ is a triangulated subcategory of $D^b(\mathcal A)$.
Let $G\in \mathcal{GP}(\mathcal A)$. Then there is a complete
projective resolution $E^\bullet$ such that ${\rm Ker}
d^{1}_{E^\bullet}\cong G$. Then the brutal truncation
$$E^\bullet_{\le 0}: = \cdots \longrightarrow E^{-2} \longrightarrow E^{-1}\stackrel
{d^{-1}_{E^\bullet}}\longrightarrow E^0 \longrightarrow 0$$ is in
$K^{-, \ b}_{\mathcal G}(\mathcal P(\mathcal A))$, and
$\rho(E^\bullet_{\le 0}) \cong {\rm Im}d^0_{E^\bullet} = {\rm Ker}
d^{1}_{E^\bullet}\cong G.$ Thus $\rho (K^{-, b}_{\mathcal
G}(\mathcal P(\mathcal A)))$ is a triangulated subcategory of
$D^b(\mathcal A)$ containing $\mathcal{GP}(\mathcal A).$ So $\langle
\mathcal{GP}(\mathcal A)\rangle\subseteq \rho (K^{-, b}_{\mathcal
G}(\mathcal P(\mathcal A))).$ Thus $\rho (K^{-, b}_{\mathcal
G}(\mathcal P(\mathcal A))) = \langle \mathcal{GP}(\mathcal
A)\rangle.$ This proves the claim.

\vskip5pt

By the claim the diagram
\[\xymatrix {K^{b}(\mathcal P(\mathcal A)) \ar[r] \ar@{=}[d] & K^{-, b}_{\mathcal G}(\mathcal P(\mathcal A))\ar[d]\ar[r]\ar[d] & K^{-,
b}(\mathcal P(\mathcal A))\ar[d] ^-{\rho}\\ K^{b}(\mathcal
P(\mathcal A)) \ar[r] & \langle \mathcal{GP}(\mathcal
A)\rangle\ar[r]&  D^b(\mathcal A)}\] commutes, where the horizontal
functors are embeddings, and the vertical functors are
triangle-equivalences.  This induces a triangle-equivalence
$$K^{-, b}(\mathcal P(\mathcal
A))/K^{-, b}_{\mathcal G}(\mathcal P(\mathcal A))\cong D^b(\mathcal
A)/\langle \mathcal{GP}(\mathcal A)\rangle = D^b_{{\rm
defect}}(\mathcal A)$$ where the last equality follows from Lemma
\ref{defect}.

\vskip5pt

$(ii)$ \ Assume that $\mathcal A$ is CM-contravariantly finite. By
Theorem \ref{DB} and Lemma \ref{verdierquotient} we get a
triangle-equivalence
$$K^{-, gpb}(\mathcal{GP}(\mathcal A))/K^b(\mathcal{GP}(\mathcal A)) \longrightarrow D^b(\mathcal
A)/ \langle \mathcal{GP}(\mathcal A)\rangle.$$ By Lemma \ref{defect}
we get a triangle-equivalence
$$D^b_{{\rm defect}}(\mathcal A) \cong K^{-, gpb}(\mathcal {GP} (\mathcal A))/K^{b}(\mathcal {GP} (\mathcal A)).$$
 This completes the proof. $\s$

\vskip5pt

\begin{rem} Let $\overline{F}: \underline {\mathcal{GP}(A)} \longrightarrow K^{-,
b}(\mathcal P(\mathcal A))/ K^{b}(\mathcal P(\mathcal A))$ be the
fully-faithful triangle functor defined in $(2.2)$. Then it is easy
to see ${\rm Im} \overline{F} = K^{-, \ b}_\mathcal G(\mathcal
P(\mathcal A))/K^b(\mathcal P(\mathcal A)).$
\end{rem}

\vskip5pt

\subsection{}  Note that a {\rm CM}-finite algebra
is {\rm CM}-contravariantly finite. By Theorem
\ref{descriptionofdeft} we get

\vskip5pt

\begin{cor} \label{dif} Let $\A$ be a {\rm CM}-finite algebra. Then there is a triangle-equivalence
$$D^b_{{\rm defect}}(\A) \cong  D_{sg}^b( {\rm \bf Aus}(\A)).$$
\end{cor}
\noindent{\bf Proof.} Put $E$ to be the direct sum of all the
pairwise non-isomorphic indecomposable Gorenstein-projective
$\A$-modules.  Then $\Hom_\A(E, -): \ \mathcal{P}( {\rm \bf
Aus}(\A)) \cong \mathcal{GP}(\A)$ as additive categories. This is
extended to a triangle-equivalences
$$K^{-, b}(\mathcal{P}({\rm \bf Aus}(\A)))/K^{b}(\mathcal{P}({\rm \bf
Aus}(\A)) \cong K^{-,
gpb}(\mathcal{GP}(\A))/K^{b}(\mathcal{GP}(\A)),$$ and hence we get
the triangle-equivalences
\begin{align*} D^b_{sg}({\rm \bf Aus}(\A)) & = D^b({\rm \bf Aus}(\A)) / K^{b}(\mathcal{P}({\rm \bf Aus}(\A))) \\ & \cong
K^{-, b}(\mathcal{P}({\rm \bf Aus}(\A)))/K^{b}(\mathcal{P}({\rm \bf
Aus}(\A)) \\ & \cong K^{-,
gpb}(\mathcal{GP}(\A))/K^{b}(\mathcal{GP}(\A))
\\& \cong D^b_{{\rm defect}}(\A),\end{align*} where the final
triangle-equivalence follows from the second triangle-equivalence in
Theorem \ref{descriptionofdeft}. This completes the proof. $\s$

\vskip5pt

\section{\bf Final Remarks} \
In fact, we have the following more general result.

\vskip5pt

\begin{thm} \label{frem}
 Suppose that $\mathcal A$ and $\mathcal A'$ are abelian categories with enough projective objects such that  $\mathcal{GP}(\mathcal A) \cong \mathcal{GP}(\mathcal A')$
as categories, and that $\mathcal A$ is {\rm CM}-contravariantly
finite. Then

\vskip5pt

$(i)$ \ There is a triangle-equivalence $$D^b(\mathcal A')\cong
K^{-, gpb}(\mathcal{GP}(\mathcal A'))/K^{b,ac}(\mathcal{GP}(\mathcal
A'));$$

\vskip5pt

$(ii)$ \ There is a triangle-equivalence
$$D^b_{defect}(\mathcal A') \cong {K^{-, gpb}(\mathcal{GP}(\mathcal A'))}/{K^{b}(\mathcal{GP}(\mathcal A'))},$$
and hence there is  a triangle-equivalence $$D^b_{defect}(\mathcal
A) \cong D^b_{defect}(\mathcal A').$$
\end{thm}

\noindent{\bf Proof.} Let $F: \mathcal{GP}(\mathcal A)  \ra
\mathcal{GP}(\mathcal A')$ be an equivalence of categories, with
quasi-inverse $F^{-1}$. Since an equivalence between additive
categories is an additive functor,  $F$ is an additive functor.

\vskip5pt

We first prove that Lemma \ref{quasi} holds also for $\mathcal A'$.
Let $P^\bullet = (P^i, d^i) \in K^{-, b}(\mathcal{P}(\mathcal A'))$,
and $N$ an integer such that ${\rm H}^{n}P^\bullet = 0$ for $n\le
N$. Consider complex $F^{-1}P^\bullet.$ Since $\mathcal{GP}(\mathcal
A)$ is  a contravariantly finite subcategory in $\mathcal A$, we can
take a right $\mathcal{GP}(\mathcal A)$-approximation $G^{N-1}
\longrightarrow {\rm Ker}F^{-1}(d^{N})$ of  ${\rm
Ker}F^{-1}(d^{N})$, and then take right $\mathcal{GP}(\mathcal
A)$-approximations  step by step (as in the proof of Lemma
\ref{quasi}) we get a complex
$$G^\bullet:  \ \cdots \longrightarrow G^{N-2} \xra{\partial^{N-2}} G^{N-1} \xra
{\partial^{N-1}} F^{-1}(P^{N}) \xra {F^{-1}(d^N)}
F^{-1}(P^{N+1})\longrightarrow \cdots.$$ By construction we have
${\rm H}^{n}\Hom_\mathcal A(G, G^\bullet) = 0$ for $G \in
\mathcal{GP}(\mathcal A)$ and $n\le N$. Since $F:
\mathcal{GP}(\mathcal A)  \ra \mathcal{GP}(\mathcal A')$ is an
equivalence, we have the complex isomorphism $\Hom_\mathcal A(G,
G^\bullet) \cong \Hom_{\mathcal A'}(FG, FG^\bullet)$ for $G \in
\mathcal{GP}(\mathcal A)$. Thus  $FG^\bullet\in K^{-,
gpb}(\mathcal{GP}(\mathcal A')).$

\vskip5pt

Since $F^{-1}$ is a  quasi-inverse of $F$,  there is a complex
$Q^\bullet$  and  a complex isomorphism
\[
\xymatrix@=13pt@R=25pt { Q^\bullet:  & \cdots  \ar[r] & F(G ^{N-2})
\ar@{=}[d]  \ar[r] &  F(G^{N-1})  \ar@{=}[d]  \ar[r]   &  P^{N}
\ar[d]  \ar[r]^{d^N}
& P^{N+1} \ar[d] \ar[r] &  \cdots\\
FG^\bullet: & \cdots \ar[r] & F(G^{N-2})  \ar[r] &  F(G^{N-1})
\ar[r] &   FF^{-1}(P^{N}) \ar[r] &  FF^{-1}(P^{N+}) \ar[r] & \cdots.
}
\]

\noindent So $Q^\bullet \in K^{-, gpb}(\mathcal{GP}(\mathcal A'))$
and ${\rm H}^{n} Q^\bullet = 0$ for  $n \le N$. Thus we have a
quasi-isomorphism
\[
\xymatrix{ P^\bullet:  & \cdots  \ar[r] & P^{N-2}  \ar[d]  \ar[r] &
P^{N-1} \ar[d]  \ar[r]   &  P^{N} \ar@{=}[d]  \ar[r]
& P^{N+1} \ar@{=}[d]  \ar[r] &  \cdots\\
Q^\bullet:  & \cdots  \ar[r] & F(G^{N-2})     \ar[r] &  F(G^{N-1})
\ar[r]   &  P^{N}   \ar[r]^{d^N} & P^{N+1}  \ar[r] &  \cdots }
\]

\noindent This proves that  Lemma \ref{quasi}$(i)$ holds for
$\mathcal A'$. Hence Lemma \ref{quasi}$(ii)$ also holds for
$\mathcal A'$.

\vskip5pt

Repeating the proof of Theorem \ref{DB}, we can prove $$
D^b(\mathcal A')\cong  K^{-, gpb}(\mathcal{GP}(\mathcal
A'))/K^{b,ac}(\mathcal{GP}(\mathcal A')).$$ This proves $(i)$.

\vskip5pt

Repeating the proof of Theorem \ref{descriptionofdeft}$(ii)$ we have
$$D^b_{defect}(\mathcal A') \cong {K^{-, gpb}(\mathcal{GP}(\mathcal
A'))}/{K^{b}(\mathcal{GP}(\mathcal A'))}.$$ Note that the
equivalence $\mathcal{GP}(\mathcal A) \cong \mathcal{GP}(\mathcal
A')$ extends to a triangle-equivalence
$${K^{-, gpb}(\mathcal{GP}(\mathcal
A))}/{K^{b}(\mathcal{GP}(\mathcal A))}\cong {K^{-,
gpb}(\mathcal{GP}(\mathcal A'))}/{K^{b}(\mathcal{GP}(\mathcal
A'))}.$$ By Theorem \ref{descriptionofdeft}$(ii)$ we have a
triangle-equivalence $D^b_{defect}(\mathcal A) \cong {K^{-,
gpb}(\mathcal{GP}(\mathcal A))}/{K^{b}(\mathcal{GP}(\mathcal A))}.$
Altogether we get a triangle-equivalence $D^b_{defect}(\mathcal A)
\cong D^b_{defect}(\mathcal A')$. This proves $(ii)$.
  $\s$

\vskip5pt

\begin{cor}
Let $\A$ and $\A'$ be Artin algebras such that $\mathcal{GP}(\A)
\cong \mathcal{GP}(\A')$ as categories.  Then $\A$ is Gorenstein if
and only if $\A'$ is Gorenstein.
\end{cor}

\noindent{\bf Proof.} If $\A$ is Gorenstein, then $\mathcal{GP}(\A)$
is contravariantly finite. By Theorem \ref{frem} we have
$$D^b_{defect}(\A') \cong  D^b_{defect}(\A) = 0,$$ thus $\A'$
is Gorenstein. $\s$

\vskip20pt

F. Kong
\par College of Math. and Statistics, \ \ Chongqing University
\par Chongqing 401331, P. R. China \ \ \ \ \
e-mail: fankong2013$\symbol{64}$gmail.com

\vskip5pt

P. Zhang\par Department of Math., \ \ Shanghai Jiao Tong University
\par Shanghai 200240, P. R. China \ \ \ \ \ \
e-mail: pzhang$\symbol{64}$sjtu.edu.cn

\medskip

\end{document}